%% file: paper.tex
\documentclass[12pt]{amsart}
\usepackage {amsmath, amssymb, amscd, graphicx, color, url,epsf}
\usepackage[margin=1in,dvips]{geometry}

\def\z{\mathbf z}
\def\Mirror{m}
\def\Os{\mathbb{O}}
\def\Xs{\mathbb{X}}
\def\orL{\vec{L}}
\def\orK{\vec{K}}
\def\tInv{\widehat\theta}
\def\lInv{\widehat\lambda}
\def\NESW{\mathcal{I}}
\def\x{\mathbf{x}}
\def\y{\mathbf{y}}
\newcommand\xUR{\mathbf{x^+}}
\newcommand\xLL{\mathbf{x^-}}
\def\R{\mathbb{R}}
\def\Z{\mathbb{Z}}
\def\tb{\textit{tb}}
\def\r{\textit{r}}
\newcommand\Rect{\mathrm{Rect}}
\newcommand\EmptyRect{\mathrm{Rect}^{\circ}}
\def\Interior{\mathrm{Int}}
\def\Torus{\mathcal{T}}
\def\Field{\mathbb{F}}
\def\sl{\textit{sl}}

\newcommand\tInvA{\widehat\theta}
\newcommand\lInvA{\widehat\lambda}
\newcommand\tInvM{\theta} 
\newcommand{\NW}{\text{\it {NW}}}
\newcommand{\NE}{\text{\it {NE}}}
\newcommand{\SW}{\text{\it {SW}}}
\newcommand{\SE}{\text{\it {SE}}}
\newcommand{\ONW}{\text{\it {O:NW}}}
\newcommand{\ONE}{\text{\it {O:NE}}}
\newcommand{\OSW}{\text{\it {O:SW}}}
\newcommand{\OSE}{\text{\it {O:SE}}}
\newcommand{\XNW}{\text{\it {X:NW}}}
\newcommand{\XNE}{\text{\it {X:NE}}}
\newcommand{\XSW}{\text{\it {X:SW}}}
\newcommand{\XSE}{\text{\it {X:SE}}}
\def\HKa{\widehat{\textit{HK}}}
\def\HKaa{\widetilde{\textit{HK}}}
\def\HKm{\textit{HK}^-}
\def\HFKa{\widehat{\textit{HFK}}}

\def\HFKm{{\textit{HFK}}^-}
\def\CFKm{{\textit{CFK}}^-}
\def\CKm{{\textit{CK}}^-}
\def\daa{\widetilde\partial}
\def\da{\widehat\partial}
\def\dm{\partial^-}
\def\daa{\widetilde\partial}
\def\Cm{C^-}
\def\Ca{\widehat{C}}
\def\CKa{\widehat{\textit{CK}}}
\def\CKaa{\widetilde{\textit{CK}}}
\newcommand\Zmod[1]{\Z/{#1}\Z}
\def\Gen{S}
\def\cm{\cdot}

\DeclareMathOperator{\Mod}{Mod}

\def\Filta{\widehat{\mathcal F}}
\def\mathcenter#1{\vcenter{\hbox{$#1$}}}

\def\mfigb#1{\mathcenter{\includegraphics[trim=-2 -2 -2 -2]{#1}}}

\input{xy}
\xyoption{all}

\newcommand\TransDiag[3]{\genfrac{}{}{0pt}{}{\mathcenter{\includegraphics[trim=-2 -2 -2 -2]{#1}}}{\hbox{\small
  $\begin{aligned}
    X &= (#2)\\
    O &= (#3)
  \end{aligned}$}}}

\newtheorem{theorem}{Theorem}
\newtheorem{proposition}{Proposition}
\newtheorem{lemma}[proposition]{Lemma}
\newtheorem{conjecture}[proposition]{Conjecture}

\begin{document}

\title[Transverse Knots Distinguished by Knot Floer Homology]{Transverse Knots Distinguished by Knot Floer Homology}

\author{Lenhard Ng}
\address{Department of Mathematics, Duke University\\ Durham, NC
27708}
\email{ng@math.duke.edu}

\author{Peter Ozsv\'ath}
\address {Department of Mathematics, Columbia University\\ New York, NY 10027}
\email {petero@math.columbia.edu}

\author{Dylan Thurston}
\address {Department of Mathematics, Barnard College, Columbia University\\ New York, NY 10027}
\email {dthurston@barnard.edu}

\begin{abstract}
  We use the recently-defined knot Floer homology invariant for transverse knots
  to show that certain pairs of transverse knots with the same self-linking number are
  not transversely isotopic. We also show that some of the algebraic refinements of
  knot Floer homology lead to refined versions of these invariants, distinguishing
  additional transversely non-isotopic knots with the same self-linking number.
\end{abstract}

\maketitle
\input{intro}
\input{comb}
\input{examples}
\input{alg}

\input{nat}

\bibliographystyle{plain}
\bibliography{biblio}

\end{document}

%% file: intro.tex
\section{Introduction}

The aim of the present paper is to study the knot Floer homology
invariant for transverse knots in $\R^3$~\cite{OST}.
Specifically, we use this invariant to distinguish transverse knots
with the same classical invariants and, indeed, we find some new
transversally nonsimple knot
types. Before stating our results, we recall some notions from
Legendrian and transverse knot theory and also knot Floer homology.

Legendrian and transverse knots play a central role in contact
geometry; see Etnyre's survey~\cite{Etnyresurvey} for further
background. For our purposes, a Legendrian knot is a knot in $\R^3$
with the property that the restriction of the standard contact form
$dz-y\,dx$ to the knot vanishes identically; a transverse knot is a
knot in $\R^3$ with the property that the restriction of the standard
contact form to the knot vanishes nowhere.

Legendrian knots in standard contact $\R^3$, modulo isotopy through
Legendrian knots, have two ``classical'' numerical invariants, the
Thurston-Bennequin number $\tb$ and the rotation number $\r$,
whereas transverse knots modulo transverse isotopy have one, the
self-linking number $\sl$. If one tries to classify Legendrian and
transverse knots in any particular topological knot type, an obvious
question arises: are the Legendrian or transverse isotopy classes
completely classified by their classical invariants? A topological
knot type is Legendrian (resp.\ transversely) simple if all
Legendrian (transverse) knots in its class are determined up to
Legendrian (transverse) isotopy by their classical invariants.

Although some knot types, including the
unknot~\cite{EliashbergFraser}, torus knots~\cite{EtnyreHondaJSG},
and the figure eight knot~\cite{EtnyreHondaJSG}, are known to be
Legendrian and transversely simple, it has been known since the work
of Chekanov and Eliashberg in the mid 1990's that not all knots are
Legendrian simple. In particular, using a invariant now called
Legendrian contact homology which counts pseudo-holomorphic curves,
Chekanov~\cite{Chekanov} produced examples of Legendrian $5_2$ knots
which have the same $\tb$ and $\r$ but are not Legendrian isotopic.
Subsequently Legendrian contact homology and other ``nonclassical''
Legendrian invariants have been used to find many other examples of
knots which are not Legendrian
simple~\cite{EpsteinFuchsMeyer,NgCLI}.

The situation with transverse knots is considerably more difficult.
Until now, the only examples of knots which are not transversely
simple have been produced using braid
theory~\cite{BirmanMenasco,MenascoMatsuda} or convex surface
theory~\cite{EtnyreHonda}. These include a family of knots of braid
index three (Birman and Menasco
\cite{BirmanMenasco,BirmanMenascoNote}) and the $(2,3)$ cable of the $(2,3)$ torus knot
(Etnyre and Honda~\cite{EtnyreHonda}).
None of these cases uses any sort of
nonclassical invariant of transverse knots.

The purpose of this paper is to show that the invariant
$\tInv(K)$ previously defined~\cite{OST} constitutes an
\emph{effective} nonclassical
invariant of transverse knots; that is, it can be used to
distinguish pairs of transverse knots with the same topological type
and self-linking number. We recall now the basic properties of $\tInv(K)$,
which takes its values in the knot Floer homology of~$K$ (in a 
sense to be made precise).

Given a knot $K\subset S^3$, knot Floer homology is a knot invariant
which is a finitely generated, bigraded Abelian group
$$\HFKa(K)=\bigoplus_{d,s}\HFKa_d(K,s),$$ whose bigraded Euler
characteristic is the symmetrized Alexander polynomial of~$K$
\cite{OS,Rasmussen}. This invariant is the homology of a chain complex
whose differentials count pseudo-holomorphic disks in a symplectic
manifold constructed from a Heegaard diagram associated to~$K$,
compare also~\cite{HolDisk}. Recent techniques have rendered the
calculations of these groups purely combinatorial
\cite{MOS,MOST,SarkarWang}.  Specifically, a grid
diagram~$G$~\cite{Brunn,Cromwell} for a knot gives rise to a bigraded
chain complex $\CKa(G)$ which can be explicitly determined from the
combinatorics of~$G$, and whose homology agrees with the knot Floer
homology groups mentioned above~\cite{MOS}.

Indeed, knot Floer homology can be
developed entirely within a combinatorial framework~\cite{MOST}. Any
two grid diagrams for a given knot can be
connected by a standard set of moves, which we call {\em grid
moves}~\cite{Cromwell} (see also Subsection~\ref{subsec:Transverse}
below). Given a sequence $M$ of grid moves from $G_1$ to $G_2$,
there is a corresponding isomorphism
$${\widehat\Phi}_M \colon \HKa(G_1) \longrightarrow \HKa(G_2).$$

A grid diagram $G$ induces a transverse realization of
its underlying knot type, and in fact any transverse knot can be
represented by a grid diagram~$G$. Moreover, there is a
restricted set of grid moves, which we call {\em transverse grid
  moves} (cf.\ Subsection~\ref{subsec:Transverse}), with the
property that two grid diagrams $G_1$ and $G_2$ represent transversely
isotopic knots if and only if $G_1$ can be connected to $G_2$ by a
sequence of transverse grid moves. This is essentially a result of
Epstein, Fuchs, and Meyer~\cite{EpsteinFuchsMeyer}, cf. also~\cite{OST}.

Combinatorial knot Floer homology and transverse knot theory meet as
follows.  The chain complex $\CKa(G)$ is equipped with a canonical
cycle $\xUR(G)$.  If $M$ is a sequence of transverse grid moves
carrying $G_1$ to $G_2$, then the induced isomorphism
${\widehat\Phi}_M$ carries the homology class of $\xUR(G_1)$ to the
homology class of $\xUR(G_2)$. Thus, the homology class $\tInv$ of
$\xUR(G)$, up to automorphisms of $\HKa(G)$, is an invariant of the
transverse isotopy class of the underlying transverse knot.

We have not yet investigated the precise dependence of
${\widehat\Phi}_M$ on the sequence of moves $M$ (though see
Section~\ref{sec:naturality}), and this may be the input required to
distinguish the three-braid examples of Birman and
Menasco~\cite{BirmanMenasco}.
(Note however that knot Floer homology
does not distinguish all of Birman and Menasco's transverse examples; in particular, all of the knots in \cite{BirmanMenascoNote}[Table III] except $11a_{240}$ have $\HFKa = 0$ and $\HFKm$ of rank $1$ in the relevant bidegree.)
But even without such an investigation, the invariant $\tInv$ can be
used to distinguish transverse isotopy classes:  for example, if
$\tInv$ vanishes for some grid representation of a transverse knot,
then it must vanish for all grid representations of transversely
isotopic knots. And indeed, we have the following:

\begin{theorem}
  \label{thm:TInvEff}
  The invariant $\tInv$ is an effective invariant of transverse
  knots.  In particular, it can be used to show that the knot types given by
  the mirrors of $10_{132}$ and $12n_{200}$ are not transversely
  simple: each of these knot types has pairs of transverse
  representatives $T_1$ and $T_2$, both with $\sl=-1$, for which
  $\tInv(T_1)=0$ and $\tInv(T_2)\neq 0$.
\end{theorem}

\noindent
This technique can also be used to distinguish transverse
representatives for the $(2,3)$ cable of the $(2,3)$ torus knot,
which was first shown to be transversely nonsimple by Etnyre and
Honda~\cite{EtnyreHonda}, see also~\cite{MenascoMatsuda}.

Some more refined invariants can also be extracted from additional
structure on knot Floer homology. Recall that knot Floer homology is
in fact the homology of the graded object associated to some
filtration of a chain complex whose homology is $\Z$, and moreover
the filtered homotopy type of this complex is a knot invariant. The
preferred isomorphisms ${\widehat\Phi}$ mentioned above are in fact
maps induced by filtered isomorphisms of the complexes.

In more concrete terms, the filtered structure immediately yields a
map $\delta_1\colon \HFKa_d(K,s)
\rightarrow\HFKa_{d-1}(K,s-1)$ which satisfies
$\delta_1^2=0$; moreover, if $${\widehat\Phi}\colon \HKa(G_1)
\longrightarrow \HKa(G_2)$$ is an isomorphism induced by grid moves,
then
$$\delta_1\circ{\widehat\Phi}={\widehat\Phi}\circ\delta_1.$$
Thus, the isomorphism class of $\delta_1\circ \tInv$ is also a transverse knot invariant.

\begin{theorem}
    \label{thm:DOneDifferential}
  The invariant $\delta_1\circ \tInv$ is an effective invariant of transverse
  knots.  In particular, it can be used to show that the pretzel knots $P(-4,-3,3)$
  and $P(-6,-3,3)$ are not transversely
  simple: each of these knot types has pairs of transverse
  representatives $T_1$ and $T_2$, both with $\sl=-1$, for which
  both $\tInv(T_1)$ and $\tInv(T_2)$ are nonzero, but
  $\delta_1\circ \tInv(T_1)=0$ while $\delta_1\circ \tInv(T_2)\neq 0$.
\end{theorem}

Knot Floer homology comes in a variety of versions. There is a version
$\HFKm$, which is the homology of a chain complex $\CFKm$ over the
ring $\Z[U]$. There is a more refined invariant $\theta^-$ of
transverse knots~\cite{OST} which is a
homology class in this variant. For our purpose, it suffices to
consider a specialization of knot Floer homology, which is a
finitely generated vector space over the field~$\Field$ with two
elements, gotten by
specializing $\CFKm$ to $\Field=\Z[U]/(2,U)$.  The
specialization to $U=0$ allows us to work with a finitely generated
Abelian group, and working in characteristic $2$ allows us to avoid sign
issues (cf.~\cite{MOST}) which demand additional computation
complexity. It seems likely that more information is contained in the
more general theory, but we do not address those issues here.
Moreover, Theorem~\ref{thm:DOneDifferential} uses $\delta_1$, but more
generally, there is an infinite sequence of maps $$\delta_k\colon
\HFKa_{d}(K,s) \longrightarrow \HFKa_{d-1}(K,s-k)$$ which could
presumably be employed to detect different transverse knot types.

The calculations underpinning Theorems~\ref{thm:TInvEff}
and~\ref{thm:DOneDifferential} above have been done by computer, using
a C program available at
\url{http://www.math.columbia.edu/~petero/transverse.html}.
Calculating knot Floer homology using the combinatorial complex is
well suited for computers; for example, Baldwin and
Gillam~\cite{BaldwinGillam} have written a program which uses this
complex to determine $\HFKa(K)$ for all knots with $11$ or fewer
crossings. Our program aims for the more modest task of determining
whether or not a given cycle is homologically trivial.  (See
Section~\ref{sec:Algorithm} for details.)  Accordingly, it is able to
handle knots of higher arc index: for example, it can be used to study
the Etnyre-Honda examples, which have grid number $17$ (the underlying
knot class has arc index $16$).

In Section~\ref{sec:Preliminaries}, we sketch the background for
this paper, starting with knot Floer homology, transverse knots, and
then the transverse invariant.  In Section~\ref{sec:Examples}, we
include the examples illustrating Theorem~\ref{thm:TInvEff},
Theorem~\ref{thm:DOneDifferential}, and the Etnyre-Honda result,
along with indications on how to find such examples.
Section~\ref{sec:Algorithm} describes the algorithm used in our
computations. In Section~\ref{sec:naturality}, we present a
conjecture on naturality in knot Floer homology and some
consequences of that conjecture, including the transverse
nonsimplicity of $7_2$ and possibly other twist knots.

\subsection*{Acknowledgements} We wish to thank John Baldwin, Marc
Culler, John Etnyre, Hiroshi Matsuda,
Olga Plamenevskaya, Jacob Rasmussen, and Zol\-t{\'a}n Szab{\'o} for
interesting conversations during the course of this work. In fact
Culler's Gridlink program~\cite{Gridlink} has proven to be an
invaluable tool for exploring examples.
LLN thanks Princeton and Columbia Universities
for their hospitality during the course of this work.
PSO was supported by NSF
grant numbers DMS-0505811 and FRG-0244663.  DPT was supported by a Sloan
Research Fellowship.



%% file: comb.tex
\section{Preliminaries}
\label{sec:Preliminaries}
\subsection{Knot Floer homology}
\label{subsec:HFK}
We review the combinatorial construction of knot Floer homology with
coefficients in $\Zmod{2}=\Field$~\cite{MOS}.

A planar \emph{grid diagram} $G$ is a diagram on an
$n \times n$ square grid in the plane, where each  square is decorated
with an $X$, an $O$, or nothing,  so that:
\begin{itemize}
\item every row contains exactly one $X$ and one $O$;
\item every column contains exactly one $X$ and one $O$.
\end {itemize}
The number $n$ is the \emph{grid number} of $G$.
Sometimes we number the $O$'s and $X$'s by
$\{O_i\}_{i=1}^n$ and $\{X_i\}_{i=1}^n$, and we denote the
two sets by~$\Os$ and~$\Xs$, respectively.

From a planar diagram, we can construct an oriented, planar link
projection by drawing horizontal segments from the $O$'s to the $X$'s
in each row, and vertical segments from the $X$'s to the $O$'s in each
column. At every crossing, the horizontal segment passes under the
vertical one. This produces a
planar diagram for an oriented link $\orL$ in $S^3$. We say that
$\orL$ has a grid presentation given by $G$.
We focus on the case where
$\orL$ is a knot~$\orK$.

If we cyclically permute the rows or columns of a grid diagram, we do
not change the knot that it represents, so we think of the
grid diagram as drawn on a torus~$\Torus$.  Let the horizontal, resp.\
vertical, \emph{(grid) circles} be the circles in between two adjacent
rows, resp.\ columns, of marked squares.  We denote the horizontal
circles by $\{\alpha_i\}_{i=1}^n$ and the vertical ones
$\{\beta_i\}_{i=1}^n$.

We associate to each toroidal
grid diagram~$G$ a chain complex
$\bigl(\CKm(G;\Field), \partial
\bigr)$ over $\Field[U_1,\dots,U_n]$.
Let $\Gen=\Gen(G)$ be the set of one-to-one correspondences between
the horizontal and vertical grid circles, which in turn can be thought
of as $n$-tuples of intersection points between the horizontal and
vertical grid circles such that no intersection point
appears on more than one horizontal or vertical grid circle. These generators
are called {\em (grid) states}.

Let $\CKm(G;\Field)$ be the free $\Field[U_1,\dots,U_n]$-module
generated by elements of $\Gen(G)$.

The complex has a bigrading, induced by two functions $A \colon \Gen
\longrightarrow \Z$ and $M \colon \Gen \longrightarrow \Z$ defined as
follows.
Given two collections $A$, $B$ of finitely many points in the plane,
let $\NESW(A,B)$ be the number of pairs $(a,b)$, where $a=(a_1,a_2)\in A$ and
$b=(b_1,b_2)\in B$ with $a_1<b_1$ and $a_2<b_2$.  Take a fundamental
domain for the torus which is cut along a horizontal and a vertical
circle, with the left and bottom edges included.  Given a generator
$\x\in\Gen$, we view $\x$ as a collection of points in this fundamental
domain.
Similarly, we view $\Os=\{O_i\}_{i=1}^n$ as a collection of
points in the plane.  Define the \textit{Maslov grading}
$$M(\x)=M_{\Os}(\x)=\NESW(\x,\x)-\NESW(\x,\Os)-\NESW(\Os,\x) +\NESW(\Os,\Os)+1.$$
Define $M_{\Xs}(\x)$ to be the same as $M_{\Os}(\x)$ with the set~$\Xs$ playing
the role of $\Os$.  We define the \textit{Alexander grading}
\[A(\x)= \frac{1}{2}\Big(M_{\Os}(\x) - M_{\Xs}(\x)\Big) - \Bigl(\frac{n-1}{2}\Bigr).\]
The module $\CKm(G;R)$ inherits a bigrading from
$M$ and $A$, with the additional convention that
multiplication by $U_i$ drops the Maslov grading by two and the
Alexander grading by one.

Given a pair of states $\x$ and $\y$, and an embedded rectangle
$r$ in $\Torus$ whose edges are arcs in the horizontal and vertical
circles, we say that $r$ connects $\x$ to $\y$ if $\x$ and $\y$ agree
along all but two horizontal circles, if all four corners of $r$ are
intersection points in $\x\cup\y$, and if the orientation induced on
each horizontal boundary component by the orientation of~$r$ inherited
from $\Torus$ goes
from a point in $\x$ to a point in $\y$.  Let $\Rect(\x,\y)$ denote
the collection of rectangles connecting $\x$ to $\y$.  If $\x,\y\in
\Gen$ agree along all but two horizontal circles, then there are exactly
two rectangles in $\Rect(\x,\y)$; otherwise $\Rect(\x,\y)=\emptyset$.
A rectangle~$r\in\Rect(\x,\y)$ is said to be \emph{empty} if
$\Interior(r)\cap\x = \emptyset$.  The space of empty rectangles
connecting~$\x$ and~$\y$ is denoted $\EmptyRect(\x,\y)$.

We endow $\CKm(G;\Field)$ with an endomorphism
$$\dm\colon
\CKm(G;\Field)\longrightarrow \CKm(G;\Field)$$ defined by
$$\dm(\x)=\sum_{\y\in\Gen}\,
\sum_{\substack{r\in\EmptyRect(\x,\y)\\X_1(r) = \dots = X_n(r) = 0}}\!\!
U_1^{O_1(r)}\cdots U_n^{O_n(r)}\cm \y,
$$
where $X_i(r)$, resp.\ $O_i(r)$, denotes the number of times $X_i$,
resp.\ $O_i$, appears in the interior of $r$.
This chain complex has two natural specializations
\begin{equation}
(\CKa(G;\Field),\da)=(\CKm(G;\Field)/(U_1=0),\dm)
\end{equation}
and
\begin{equation}
(\CKaa(G;\Field),\daa)=(\CKm(G;\Field)/(U_1=\dots=U_n=0),\dm).
\label{def:CKaa}
\end{equation}

Knot Floer homology comes in two natural forms, $\HFKa(K)$ and
$\HFKm(K)$, the first of which is a vector space over $\Field$, and
the second of which is a module over $\Field[U]$.
Let $V$ be the two-dimensional bigraded vector space spanned
by one generator in bigrading $(-1,-1)$ and another in bigrading $(0,0)$.
Then the two forms of knot Floer homology are related to the above
specializations, according to the
following special case of a more general result:
\begin{theorem}[\cite{MOS}]
  \label{thm:KnotFloerHomology} Fix a grid presentation $G$ of a
  knot $K$, with grid number $n$. Then the homology
  groups $H_*(\CKm(G),\dm)$ and $H_*(\CKa(G),\da)$
  are the knot invariants $\HFKm(K)$ and $\HFKa(K)$ respectively.
  The homology $H_*(\CKaa(G),\daa)$ is isomorphic to
  $\HFKa(K)\otimes V^{\otimes (n-1)}$.
\end{theorem}

There are refinements of the above construction; we describe one
that is useful here.  Consider the complex $\Cm(G;\Field)$ over
$\Field[U_1,\dots,U_n]$, whose underlying module agrees with
$\CKm(G;\Field)$, but is equipped with the endomorphism
$$\dm =
\sum_{k=0}^{\infty} \partial_k^-,$$
where here
$$\partial_k^-\colon \CKm_{d}(K,s)\longrightarrow \CKm_{d-1}(K,s-k)$$
is given by
\[\partial_k^-(\x)=\sum_{\y\in\Gen}\,
\sum_{\substack{r\in\EmptyRect(\x,\y)\\
\sum X_i(r) = k
}}\!
U_1^{O_1(r)}\cdots U_n^{O_n(r)}\cdot \y.
\]
The map $\dm$ is a differential on the complex $\Cm(G;\Field)$
equipped with its Maslov grading (i.e., it decreases the Maslov grading
by one), which respects the filtration on $\Cm(G;\Field)$ induced by
$A$ (the {\em Alexander filtration}).  Similarly define $\da_k$ and $\daa_k$.

The main result of~\cite{MOS} identifies the filtered
quasi-isomorphism type of $(\Cm,\dm)$ with the topologically
invariant ``knot filtration'' of~\cite{OS,Rasmussen}. Concretely, if
grid diagrams $G_1$ and $G_2$ represent the same knot, then in fact
the filtered complexes $\Cm(G_1,\dm)$ and $\Cm(G_2,\dm)$ are
filtered quasi-isomorphic.

$\CKm(K)$ is the associated graded
object for the Alexander filtration of $\Cm(K)$, and so its
homology is a knot invariant, as stated in
Theorem~\ref{thm:KnotFloerHomology}.  But the filtered
quasi-isomorphism type of a complex has other invariants: indeed, the
entire Leray spectral sequence is preserved under filtered
quasi-isomorphisms (cf.~\cite{UsersGuide}). In our case, this principle
can be formulated as follows:

\begin{proposition}
  \label{prop:SpectralSequences}
  Given a grid diagram $G$, inductively
  define chain complexes $(E_k(G),\delta_k)$ by
  \begin{align*}
    (E_0^-(G),\delta_0^-)&=(\CKm(G),\dm) \\
    (E_k^-(G),\delta_k^-)&=(H_*(E_{k-1}^-(G),\delta_{k-1}^-),[\dm_k]).
  \end{align*}
  If $G_1$ and $G_2$ represent isotopic knots, then there are
  isomorphisms
  $$\Phi_k^-\colon (E_k^-(G_1),\delta_k)\longrightarrow (E_k^-(G_2),\delta_k)$$
  for all $k\geq 0$.
  We can define a similar spectral sequence for $\CKa(G)$: given
  a grid diagram $G$,
  inductively define chain complexes $({\widehat
    E}_k(G),{\widehat\delta}_k)$ by
  \begin{align*}
    ({\widehat E}_0(G),{\widehat \delta})&=(\CKa(G),\widehat \partial) \\
    ({\widehat E}_k(G),{\widehat\delta}_k)&=
    (H_*({\widehat E}_{k-1}(G),{\widehat\delta}_{k-1}),[{\widehat \partial}_k]).
  \end{align*}
  If $G_1$ and $G_2$ represent isotopic knots, then there are
  isomorphisms of chain complexes
  $${\widehat \Phi}_k\colon ({\widehat E}_k(G_1),{\widehat\delta}_k)
  \longrightarrow ({\widehat E}_k(G_2),{\widehat \delta}_k)$$
  for all $k\geq 0$.
  Furthermore, the canonical map $i\colon \CKm(G)\longrightarrow \CKa(G)$
  induces a map of spectral sequences making the following diagram commute:
  \[
  \begin{CD}
    (E_k^-(G_1),\delta^-_k) @>{\Phi_k^-}>> (E_k^-(G_2),\delta^-_k) \\
    @V{i_k}VV @VV{i_k}V \\
    ({\widehat E}_k(G_1),{\widehat\delta}_k) @>{{\widehat \Phi}_k}>>
    ({\widehat E}_k(G_2),{\widehat\delta}_k).\\
  \end{CD}
  \]
  Finally, we can define a spectral sequence
  \begin{align*}
    ({\widetilde E}_0(G),{\widetilde \delta})&=(\CKaa(G),\widetilde\partial) \\
    ({\widetilde E}_k(G),{\widetilde\delta}_k)&=
    (H_*({\widetilde E}_{k-1}(G),{\widetilde\delta}_{k-1}),
    [{\widetilde\partial}_k]).
  \end{align*}
  In this case, the canonical map $j\colon \CKa(G)\longrightarrow \CKaa(G)$
  induces injective chain maps for all $k\geq 0$:
  $$j_k\colon ({\widehat E}_k(G),{\widehat \delta}) \longrightarrow
  ({\widetilde E}_k(G),{\widetilde \delta})
  \cong ({\widehat E}_k(G),{\widehat \delta})\otimes V^{\otimes (n-1)}.$$
\end{proposition}

\begin{proof}
  The filtered quasi-isomorphism class
  of the module $\Cm(G;\Field)$ is a knot invariant~\cite{MOS}.
The first two spectral sequences are naturally associated
  to this quasi-isomorphism type. Properties of the third spectral sequence also follow from general
  algebraic principles; cf.\ \cite[Lemma~\ref{MOST:lemma:UniversalCoefficients}]{MOST}.
\end{proof}

\subsection{Transverse knots and grid diagrams}
\label{subsec:Transverse}

We now review the relation between grid diagrams and Legendrian and
transverse knots.
According to Cromwell~\cite{Cromwell} (see also~\cite{Dynnikov}),
two grid diagrams $G_1$ and $G_2$ on the torus represent isotopic knots in $S^3$
if and only if they can be connected by a finite sequence of the following
{\em grid moves}:
\begin{description}
\item[Commutation]  For any pair of consecutive columns of~$G$
  so that the~$X$ and~$O$ from one column do not separate the~$X$
  and~$O$ from the other column, switch the decorations of these two
  columns.  There is also a
  similar move using rows rather than columns.
\item[Destabilization] For a corner $c$ which is shared by a pair of
  vertically-stacked squares marked with an $X$ and $O$, we delete the
  horizontal and vertical circles containing $c$, and remove the
  markings of the initial $X$ and $O$ (both of which mark now the same
  square in the destabilized diagram). We further assume
  that one of the initial squares $X$ and $O$ meets an additional square
  marked by an $X$ or an $O$. 
\item[Stabilization] The inverse of destabilization.
\end{description}
In~\cite{MOST}, an independent proof of the topological invariance of
knot Floer homology is provided, by exhibiting explicit filtered
quasi-isomorphisms between the complexes $\Cm(G)$, as the grid
undergoes each of the above grid moves.

It will be convenient to classify (de)stabilization moves according to
the local configuration of $X$'s and $O$'s. For any destabilization,
there are three marked squares in the original diagram sharing one
corner. There are two pieces of data to keep track of: the marking
shared by two of these three squares (i.e., an $X$ or an $O$), and the
placement of the \emph{unmarked} square relative to the shared
corner, either $\NW$, $\SW$, $\SE$, or $\NE$.
Stabilizations then fall into eight types. Of these, the types $\OSE$,
$\ONE$, $\ONW$, $\OSW$ are equivalent modulo commutation moves to
a stabilization of type $\XNW$, $\XSW$, $\XSE$, $\XNE$, respectively
\cite[Lemma~\ref{LegInv:cor:Stabilization}]{OST}.

There are restricted sets of moves to describe Legendrian and
transverse knots. Before describing these, we make a quick digression into
Legendrian knots.

Recall that Legendrian knots in $\mathbb{R}^3$, endowed with the
standard contact form $dz-y\,dx$, are knots along which $dz-y\,dx$
vanishes identically. These can be studied via their front
projections: images under the map $(x,y,z)\mapsto (x,z)$. A
generic Legendrian front projection is a curve in the $(x,z)$ plane
which has no vertical tangencies and is smooth away from finitely
many cusps and double-point crossings.

As explained earlier, a planar grid diagram $G$ induces a projection
for a knot $K$. It also induces a Legendrian front projection for the
mirror~$\Mirror(K)$ of~$K$: Starting with
the projection of $K$ corresponding to $G$, reverse all crossings (so
that horizontal segments cross over vertical ones), smooth all
northwest and southeast corners, turn southwest and northeast corners
into cusps, and tilt the diagram $45^\circ$ clockwise (so that the NE,
resp.\ SW, corners become right, resp.\ left, cusps).  With the $x$
axis as horizontal and $z$ as vertical, this gives a Legendrian front
projection for the mirror of the knot~$K$ described by~$G$.
Conversely, any Legendrian knot is Legendrian isotopic to the front
obtained from some grid diagram.

\begin{figure}
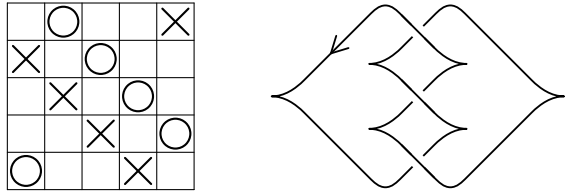

\includegraphics[height=1in]{draws/knots.60}
\qquad
\includegraphics[height=1in]{draws/knots.61}
\caption{
Grid diagram for a right-handed trefoil and the corresponding oriented
Legendrian left-handed trefoil knot.
}
\label{trefoil}
\end{figure}

There are two classical invariants of oriented Legendrian knots
modulo isotopy through Legendrian knots: Thurston-Bennequin number
$\tb$ and rotation number $\r$. In terms of grid diagrams, these are
defined as follows. Consider a grid diagram $G$ corresponding to an
oriented Legendrian knot $L$. Let $\operatorname{wr}(G)$ be the
writhe of the knot projection given by $G$, i.e., the number of
positive crossings minus the number of negative crossings; note that
because of crossing changes, this is negative the writhe of the
front of $m(L)$ considered as a knot diagram. Let $\# \NE(G)$ denote
the number of northeast corners in the knot projection given by $G$,
let $\# \NE_X(G)$ be the number of these corners occupied by $X$'s,
and similarly define $\# \NE_O(G)$, $\# \SW_X(G)$, and $\#
\SW_O(G)$. Then
\begin{align*}
\tb(L) &= -\operatorname{wr}(G) - \# \NE(G) \\
\r(L) &= \frac{1}{2} \bigl(
\# \NE_X(G) - \# \NE_O(G) - \# \SW_X(G) + \#
\SW_O(G) \bigr).
\end{align*}

For example, Figure~\ref{trefoil} gives a grid diagram for the
right-handed trefoil, yielding a
Legendrian left-handed trefoil. This
diagram has $\operatorname{wr} = 3$, $\# \NE = 3$, $\#
\NE_X = 1$, $\# \NE_O = 2$, $\# \SW_X = 2$,
and $\# \SW_O = 1$, and thus $\tb = -6$ and $\r = -1$. In the
sequel we will suppress crossing information in Legendrian front
diagrams; at any crossing, the negatively sloped strand passes over
the positively sloped strand.

For grid diagrams, Legendrian isotopy may be expressed as follows: two
grid diagrams correspond to Legendrian isotopic knots if and only if
they can be related by commutation moves and de/stabilization moves of
type $\XNW$ and $\XSE$. $\OSE$ and $\ONW$ moves can also be included
here if desired. The other stabilization moves, while not changing
topological knot type, do change Legendrian isotopy class: $\XNE$ (or
$\OSW$) is called \emph{positive stabilization}, while $\XSW$ (or
$\ONE$) is called \emph{negative stabilization}. These
stabilizations both decrease $\tb$ by $1$; positive stabilization
increases $\r$ by $1$, while negative stabilization decreases $\r$ by
$1$. In terms of front
projections, stabilization replaces a smooth section of the front by a
zigzag, situated to change $\r$ by $\pm 1$.

Closely related to Legendrian knots are transverse knots, which are
everywhere transverse to the contact $2$-plane field $\ker
(dz-y\,dx)$. Each transverse knot inherits a natural orientation
from the coorientation of the contact structure.  That is, the
contact $1$-form $dz-y\,dx$ evaluated on tangent vectors to the knot
is always positive.
Any oriented Legendrian knot $L$ can be perturbed in the
$C^\infty$ topology to a transverse knot
by pushing it along its length an arbitrary small amount in a generic
direction transverse to the
contact planes. The resulting transverse knot is called the
\emph{positive} or \emph{negative transverse pushoff} of $L$,
written $L^+$ or $L^-$, depending on whether its orientation agrees or
disagrees with the orientation from $L$. Pushing $L$ in opposite
directions yields transverse pushoffs of opposite sign, and $L^-$ is
the positive transverse pushoff of the orientation reverse of $L$.

Both $L^+$ and $L^-$ are well-defined up to isotopy through
transverse knots, and Legendrian isotopic knots have pushoffs which
are transversely isotopic. Indeed, the set of transverse knots up to
transverse isotopy is naturally identified, through the
correspondence $L^+ \leftrightarrow L$, with the set of Legendrian
knots up to Legendrian isotopy and negative stabilization
\cite{EpsteinFuchsMeyer}.
This fact readily leads to the following characterization of
transverse isotopy in terms of grid
diagrams~\cite[Corollary~\ref{LegInv:cor:TransverseLinks}]{OST}:



\begin{proposition}
  Two grid diagrams represent Legendrian links whose positive
  transverse pushoffs are transversely isotopic if and only if they
  can be connected by a sequence of
  commutation and de/stabilization moves of types
  $\XNW$, $\XSE$, and $\XSW$.
\end{proposition}

Negative stabilization does not change $\tb-\r$. We thus define
the \emph{self-linking number} of the
corresponding positive pushoff transverse knot:
\[
\sl(L^+) = \tb(L) - \r(L).
\]
The self-linking number is then invariant under transverse isotopy.
We also have $\sl(L^-) = \tb(L) + \r(L)$.


\subsection{The transverse invariant}

We briefly recall now the transverse invariant~\cite{OST}.

Given a grid diagram $G$, consider the chain $\x^+(G)$ which occupies the
upper right corner of each square marked with an $X$. This element $\x^+(G)$
is easily seen to be a cycle in $\CKm(G)$. In fact,
the homology class of $\x^+(G)$ is an invariant of the
transverse isotopy class of the transverse knot described by $G$~\cite{OST}.
This is proved by showing that if $G_1$ and $G_2$ are grid diagrams
which differ
by a transverse grid move, then the corresponding (filtered) quasi-isomorphism
from $\CKm(G_1)$ to $\CKm(G_2)$ identified in~\cite{MOST} carries the homology
class of $\x^+(G_1)$ to that of $\x^+(G_2)$. Thus, if $G$ represents
the topological knot $K$ and the Legendrian knot $L$ of type $m(K)$,
then the underlying homology
class, $\lambda_+(L)\in \HFKm(K)$, is an invariant of the transverse
isotopy class of $L^+$. This results in a transverse invariant
$\tInvM$ defined by $\tInvM(L^+) = \lambda_+(L)$.

We can alternatively consider the chain $\x^-(G)$ which occupies the
lower left corner of each square marked with an $X$. The homology class of
this element $\lambda_-(L)$ is then equal to $\tInvM(L^-)$.

\begin{theorem}[\cite{OST}]
  \label{thm:OST}
  The homology classes $\lambda_\pm(L)$ are supported in
  $\HKm_{d}(\Mirror(L),s)$, where $d=\sl(L^\pm)+1$, $2s=d$.  Moreover,
  if $L_1^+$
  and $L_2^+$ represent transversely isotopic transverse knots, then
  there is a filtered quasi-isomorphism $\Phi^-\colon
  \Cm(\Mirror(L_1^+))\longrightarrow\Cm(\Mirror(L_2^+))$ whose induced map on homology
  $$\phi^-\colon \HKm(\Mirror(L_1^+))\longrightarrow\HKm(\Mirror(L_2^+))$$ carries
  $\lambda_+(L_1^+)$ to $\lambda_+(L_2^+)$. An analogous result holds
  for $\lambda_-$ if $L_1^-$ and $L_2^-$ are transversely isotopic.
\end{theorem}

In practice, it is more convenient to work with $\HKa$ rather group
$\HKm$, which is infinitely generated over~$\Field$.
Correspondingly, we let $\lInvA_+$, $\lInvA_-$, $\tInvA$ denote the
images of $\lambda_+$, $\lambda_-$, $\tInvM$, respectively, under
the natural map
$$i\colon \HKm(\Mirror(L)) \longrightarrow \HKa(\Mirror(L)).$$
These images are also invariants of the respective Legendrian and
transverse knot types, in view of Theorem~\ref{thm:OST} and
Proposition~\ref{prop:SpectralSequences}.


%% file: examples.tex
\section{Examples}
\label{sec:Examples}


In this section, we prove Theorems~\ref{thm:TInvEff}
and~\ref{thm:DOneDifferential}, showing that $\tInv$
and $\delta_1\circ\tInv$ are effective nonclassical transverse
invariants. The proofs consist of presenting pairs of Legendrian
knots in the relevant knot types and appealing to the computer
program. We also include some remarks for each example, and a
concluding subsection explaining the strategy used to find our examples.

For ease of reference, we collect our conventions here: $L^{\pm}$
are the positive and negative transverse pushoffs of Legendrian
$L$; $\sl(L^{\pm}) = \tb(L) \mp \r(L)$, $\theta(L^{\pm}) =
\lambda_{\pm}(L)$, and Legendrian knots are negatively/positively
stably isotopic if and only if their positive/negative transverse
pushoffs are isotopic.

\subsection{$m(10_{132})$ and $m(12n_{200})$}
\label{ssec:10132}

\begin{figure}
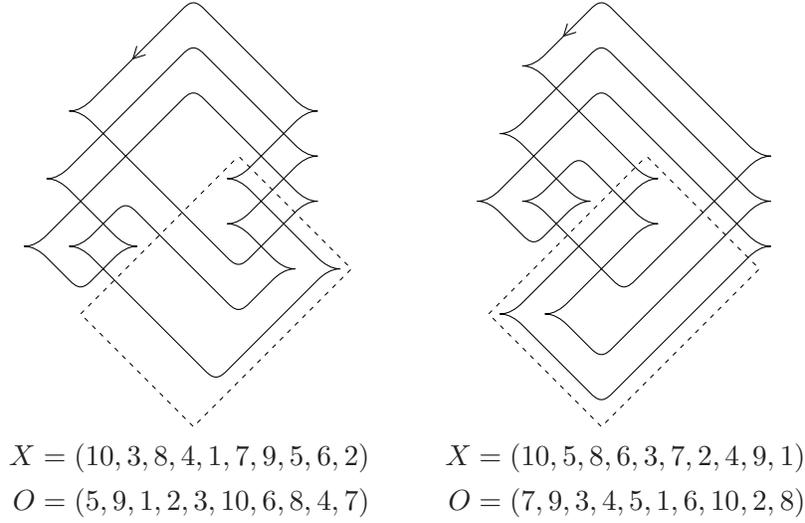

\[
\TransDiag{draws/knots.11}
  {10,3,8,4,1,7,9,5,6,2}
  {5,9,1,2,3,10,6,8,4,7}
\qquad
\TransDiag{draws/knots.16}
  {10,5,8,6,3,7,2,4,9,1}
  {7,9,3,4,5,1,6,10,2,8}
\]
\caption{
Legendrian $m(10_{132})$ knots $L_1$ (left) and $L_2$ (right). As
pictured, these knots differ only in the dashed boxes, which are
each topologically a negative half-twist on three strands.
}
\label{10132}
\end{figure}

Let $L_1$ and $L_2$ denote the oriented Legendrian knots of
topological type $m(10_{132})$ (the mirror of $10_{132}$) whose
front projections are given in Figure~\ref{10132}.\footnote{$L_1$
is also a Legendrian representative of
$m(10_{132})$, the only topological knot with $10$ or fewer
crossings for which the maximal value of $\tb$ is currently
unknown~\cite{NgKhovanov}.} Both $L_1$ and $L_2$ have $\tb=-1$ and $\r=0$, and hence
the transverse pushoffs $L_1^\pm$ and $L_2^\pm$ all have $\sl=-1$.
Note that $L_1$ and $L_2$ differ only within the dashed boxes.

The depictions of $L_1$ and $L_2$ in Figure~\ref{10132} have been
chosen to be easy to translate to grid diagrams. We can represent
an $n\times n$ grid diagram by two $n$-tuples $X = (x_1,\dots,x_n)$
and $O = (o_1,\dots,o_n)$, both permutations of $(1,\dots,n)$, so that
column $i$ contains an $X$ in row $x_i$ and an $O$ in row $o_i$, where
we number rows from bottom to top and columns from left to
right. Then the tuples for $L_1$ and $L_2$ are
given below the respective diagrams.
(It is possible to represent $L_1$ and $L_2$ by diagrams with grid
number $9$ rather than $10$, but this obscures their similarity.)

We remark that $L_2$ is Legendrian isotopic to the orientation reverse
$-L_1$ of $L_1$, and thus $L_2^\pm$ is transversely isotopic to
$L_1^\mp$. This can be seen as follows. Reflection of an oriented front
in the vertical axis does not change Legendrian isotopy class, since
it corresponds to rotating the $xy$ projection of the knot by
$180^\circ$. Reflecting $L_2$ in the vertical axis, and then
pushing the pattern in the dashed box around the knot as in the
Legendrian satellite construction \cite{NgTraynor}, yields $-L_1$.

\begin{figure}
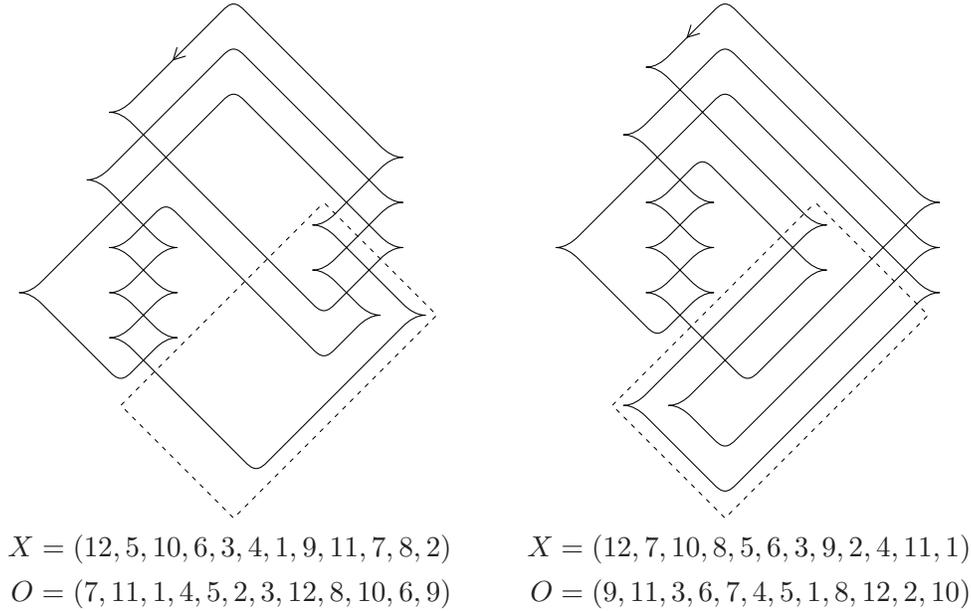

\[
\TransDiag{draws/knots.50}
  {12,5,10,6,3,4,1,9,11,7,8,2}
  {7,11,1,4,5,2,3,12,8,10,6,9}
\qquad
\TransDiag{draws/knots.51}
  {12,7,10,8,5,6,3,9,2,4,11,1}
  {9,11,3,6,7,4,5,1,8,12,2,10}
\]
\caption{
Legendrian $m(12n_{200})$ knots $L_1'$ (left) and $L_2'$ (right).
}
\label{12n200}
\end{figure}

In a similar vein, let $L_1'$ and $L_2'$ denote the Legendrian
$m(12n_{200})$ knots depicted in Figure~\ref{12n200}, also with
$\tb=-1$ and $r=0$. Then $(L_1')^\pm$
and $(L_2')^\pm$ all have $\sl=-1$, $L_2'$ is Legendrian isotopic to
$-L_1'$, and $(L_2')^\pm$ is transversely isotopic to $(L_1')^\mp$.

\begin{proposition}
$L_1^+$ and $L_2^+=L_1^-$ are not transversely isotopic; $(L_1')^+$ and
$(L_2')^+=(L_1')^-$ are not transversely isotopic.
\label{prop:10132}
\end{proposition}

\begin{proof}
The computer program tells us that $\tInv(L_1^+)$ is
null-homologous in $\HFKa_0(m(10_{132}),0)$ while $\tInv(L_2^+)$ is
not null-homologous; similarly, $\tInv((L_1')^+) =
0$ in $\HFKa_0(m(12n_{200}),0)$ while $\tInv((L_2')^+) \neq 0$.
\end{proof}

We remark that the Plamenevskaya invariant of transverse knots from
Khovanov homology \cite{Plamenevskaya} vanishes for
$L_1^{\pm}=L_2^{\mp}$ and $(L_1')^{\pm} = (L_2')^{\mp}$, because
in each case it lies in a trivial graded
summand of Khovanov homology.

Presumably $L_1$ and $L_2$ are part of a family of knots, obtained
by adding further full twists, which are not transversely simple,
but establishing this result for the entire family might be
computationally difficult.


A corollary of Proposition~\ref{prop:10132} is that $L_1$ and $-L_1$
are not Legendrian isotopic. We note that Legendrian contact
homology, which in general is quite good at distinguishing different
Legendrian knots, does not distinguish between $L_1$ and $-L_1$; the
$\Z[t^{\pm 1}]$ differential graded algebras \cite{ENS} for $L_1$
and $-L_1$, which in theory could be distinct, are easily shown to
be stable tame isomorphic. This observation also holds for $L_1'$, as
well as the knots $L_1$ and $L_2'$ defined in the next subsection.

\subsection{$P(-4,-3,3)$ and $P(-6,-3,3)$}
\label{ssec:pretzel}

In the previous section, we presented knots which were distinguished
straightaway by $\tInv$. Here we give examples of transverse
pretzel knots which are not distinguished by $\tInv$, but
rather by $\delta_1 \circ \tInv$.

\begin{figure}
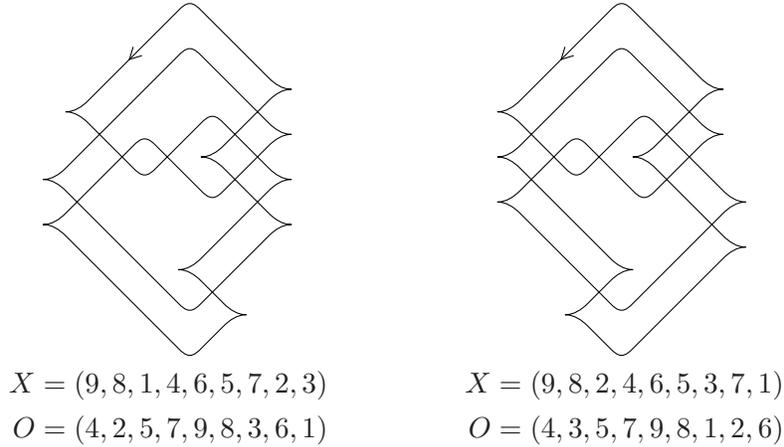

\[
\TransDiag{draws/knots.30}
  {9,8,1,4,6,5,7,2,3}
  {4,2,5,7,9,8,3,6,1}
\qquad\qquad
\TransDiag{draws/knots.32}
  {9,8,2,4,6,5,3,7,1}
  {4,3,5,7,9,8,1,2,6}
\]
\caption{Legendrian $P(-4,-3,3)$ pretzel knots $L_1$ (left) and
$L_2$ (right).}
\label{pretzel433}
\end{figure}

Let $L_1$ and $L_2$ be the Legendrian $P(-4,-3,3)=m(10_{140})$
pretzel knots shown in Figure~\ref{pretzel433}, and let $L_1'$ and
$L_2'$ be the Legendrian $P(-6,-3,3)=12n_{582}$ pretzel knots shown
in Figure~\ref{pretzel633}. All of these knots have $\tb=-1$ and
$\r=0$, and thus $L_1^{\pm}$, $L_2^{\pm}$, $(L_1')^{\pm}$, and
$(L_2')^{\pm}$ all have $\sl = -1$. Note that $L_1$ and $L_2$
(resp. $L_1'$ and $L_2'$) differ by the placement of a half
twist, \`a la the ``Eliashberg knots'' of \cite{EpsteinFuchsMeyer},
and it is easy to check that $L_1$ and $L_2$ (resp.\ $L_1'$ and
$L_2'$) are Legendrian isotopic after one positive stabilization.
However, the following result shows that they are not Legendrian
isotopic after any number of negative stabilizations.

\begin{figure}
\[
\TransDiag{draws/knots.35}
   {11,10, 4, 5, 1, 6, 8, 7, 9, 2, 3}
   { 6, 5, 7, 2, 4, 9,11,10, 3, 8, 1}
\qquad\qquad
\TransDiag{draws/knots.37}
   {11,10, 4, 5, 2, 6, 8, 7, 3, 9, 1}
   { 6, 5, 7, 3, 4, 9,11,10, 1, 2, 8}
\]
\caption{
Legendrian $P(-6,-3,3)$ pretzel knots $L_1'$ (left) and $L_2'$
(right).
}
\label{pretzel633}
\end{figure}

\begin{proposition}
Each of the following pairs of knots is not transversely isotopic: $L_1^+$
and $L_2^+$; $(L_1')^+$ and $(L_2')^+$; $L_1^+$ and $L_1^-$;
$(L_2')^+$ and $(L_2')^-$.
\label{pretzel}
\end{proposition}

\begin{proof}
The computer program finds that in
$\HFKa_{-1}(P(-4,-3,3),-1)$,
\begin{align*}
\delta_1(\tInv(L_1^+)) &= \delta_1(\lInv_+(L_1)) = 0\\
\delta_1(\tInv(L_2^+)) &= \delta_1(\lInv_+(L_2)) \neq 0\\
\delta_1(\tInv(L_1^-)) &= \delta_1(\lInv_-(L_1)) \neq 0;
\end{align*}
similarly in $\HFKa_{-1}(P(-6,-3,3),-1)$,
\begin{align*}
\delta_1(\tInv((L_1')^+)) &= \delta_1(\tInv((L_2')^-)=0\\
\delta_1(\tInv((L_2')^+)) &\neq 0.
\end{align*}
 (On the other hand, $\tInv(L_1^\pm)$, $\tInv(L_2^\pm)$,
$\tInv((L_1')^\pm)$, $\tInv((L_2')^\pm)$ are all nonzero in
$\HFKa$.)
\end{proof}

One might guess that there is a generalization of
Proposition~\ref{pretzel} to show that no pretzel knot
$P(-2n,-3,3)$ for $n\geq 2$ is transversely simple.

As in the previous subsection, the Plamenevskaya transverse invariant
\cite{Plamenevskaya} does not distinguish any of the examples in
this subsection, though for a different reason. Here the Khovanov homology in the relevant bigrading has rank $1$, and all of the
transverse knots involved have
quasipositive braid representatives. To obtain a braid from a grid diagram, connect $X$'s and $O$'s as usual, but for each vertical segment with $X$ above $O$, replace the segment by a vertical segment that starts at the $X$, goes up through the top of the diagram, and then returns to the $O$ through the bottom of the diagram. (All horizontal segments lie over all vertical segments.) We obtain a braid going from bottom to top, and it can be proven that the transverse knots represented by the grid diagram and by the corresponding braid are the same. For instance, $L_1^+$ and $L_2^+$ are represented by the quasipositive braids
$\sigma_3^{-1}\sigma_2\sigma_3^2\sigma_1^2\sigma_2^{-1}\sigma_1\sigma_2\sigma_1^{-2}$,
and 
$\sigma_3 \sigma_2 \sigma_1^2 \sigma_3^{-1} \sigma_2^{-1} \sigma_1 \sigma_2 \sigma_3 \sigma_1^{-2}$, respectively.

\subsection{The Etnyre-Honda cable example}

In \cite{EtnyreHonda}, Etnyre and Honda describe a
knot class, the $(2,3)$ cable of the $(2,3)$ torus knot, which has
a Legendrian representative which does not maximize
Thurston-Bennequin number but is nonetheless not
destabilizable. Their classification of Legendrian knots in this knot
class includes the following result.

\begin{proposition}[Etnyre-Honda \cite{EtnyreHonda}]
There are nonisotopic Legendrian knots $L_1$ and $L_2$,
topologically the $(2,3)$
cable of the $(2,3)$ torus knot, both with $\tb=5$ and $\r=2$, for
which $L_1$ is the positive stabilization of a Legendrian knot while
$L_2$ is not. Furthermore, $L_1$ and $L_2$ are not Legendrian
isotopic after any number of negative stabilizations:
$L_1^+$ and $L_2^+$ are not transversely isotopic.
\label{EH}
\end{proposition}

\begin{figure}
\[
\TransDiag{draws/knots.20}
  {10,17,9,14,6,12,11,1,7,15,13,3,2,8,16,5,4}
  {1,14,15,7,13,3,2,8,16,5,4,10,9,17,6,12,11}
\]
\[
\TransDiag{draws/knots.22}
  {10,9,16,17,6,12,11,1,7,14,13,3,2,8,15,5,4}
  {1,17,7,14,13,3,2,8,15,5,4,10,9,16,6,12,11}
\]
\caption{Legendrian fronts for $L_1$ (top) and $L_2$
  (bottom), which are both $(2,3)$ cables of the $(2,3)$ torus knot.
  These examples are derived from diagrams of Menasco and Matsuda
  \cite{MenascoMatsuda}, Figures 16 and 15, respectively. Note that
  $L_1$ and $L_2$ differ only in the indicated regions.}
\label{menasco}
\end{figure}

Menasco and Matsuda~\cite{MenascoMatsuda} have presented explicit
forms for $L_1$ and $L_2$; equivalent but slightly modified
versions, arranged to emphasize the local change that relates them,
are given in Figure~\ref{menasco}.  Since $L_1$ is a positive
stabilization (as readily checked by Gridlink), it follows from
\cite[Theorem~1.3]{OST} that $\tInv(L_1^+) = 0$ in
$\HFKa$; this can be confirmed by computer. On the other
hand, the computer program verifies that the image of
$\tInv(L_2^+)$ in $\HFKa$ is nonzero. Thus one can use $\tInv$
to reprove Proposition~\ref{EH}.

\subsection{Finding transversely nonsimple knots}

To conclude this section, we give a heuristic explanation for how
the transverse examples in Subsections~\ref{ssec:10132}
and~\ref{ssec:pretzel} were found. The techniques described here
should allow the interested reader to find other examples of knot
types which are transversely nonsimple.

We first describe how to find knot types for which the nonvanishing
of $\tInv$ might be applied to distinguish transverse
representatives. The key observation here is that such knot types must
be {\em thick}; that is, their knot Floer homology $\HFKa$ must be
supported on more than one diagonal. Indeed, we have the following
result.

\begin{proposition}
Let $T$ be a transverse knot in a thin knot type $K$. Then
$\tInvA(T) \neq 0$ if and only if $\sl(T) = 2\tau(K)-1$.
\label{prop:tau}
\end{proposition}

\noindent Note that by Plamenevskaya \cite{Plam:HFK}, any transverse
knot $T$ of type $K$ satisfies $\sl(T) \leq 2\tau(K)-1$; thus
Proposition~\ref{prop:tau} states that the $\tInvA$ invariant of a
transverse knot of thin type is nonzero if and only if the $\tau$
bound is sharp.

\begin{proof}[Proof of Proposition~\ref{prop:tau}]
Let $A$ and $M$ denote Alexander and Maslov grading, respectively,
in both $\HFKa$ and $\HFKm$. If $K$ is thin, it is an easy exercise in
homological algebra to see that, along the line $M=2A$, $\HFKm(m(K))$
has the following form: it consists of the direct sum of a free module 
over $\Field[U]$ generated by one generator in bidegree 
$(A,M) = (\tau(K),2\tau(K))$,
and a summand in bidegree $(\tau(K),2\tau(K))$ which is
annihilated by multiplication by $U$.
If $T$ is a transverse knot of type $K$, then
$\tInvM(T)$ lies in bidegree $((\sl(T)+1)/2,\sl(T)+1)$ 
and is non-$U$-torsion
\cite[Theorem~\ref{LegInv:thm:NonTorsion}]{OST}; 
it follows that $\tInvM(T)$ is in the image of multiplication by $U$ 
if and only if
$\sl(T) < 2\tau(K)-1$. Since $\tInvA(T) = 0$ if and only if
$\tInvM(T)$ is in the image of $U$, the proposition follows.
\end{proof}

Now suppose that there are transverse knots $T_1$, $T_2$ in type $K$
with the same $\sl$, for which $\tInvA(T_1) = 0$ and $\tInvA(T_2)
\neq 0$. By Proposition~\ref{prop:tau}, $K$ must be thick. In
addition, $\tInvM(T_1)$ and $\tInvM(T_2)$ are nonzero, unequal
elements in $\HFKm(m(K))$ in grading $((\sl+1)/2,\sl+1)$. In
particular, there must be a group $\HFKm_d(m(K),s)$ with $d=2s$ which
has rank at least $2$ over $\Field$.

\begin{figure}
\[
\xymatrix@=25pt@!0{ &&&M&&&
&& &&&M&&& \\
&&&&1 \ar[dl]_1 & 1 \ar[dl]^1 &
&& &&&&1 \ar@{.>}[ddl] & 1 &\\
&&& 1 \ar[uu] \ar@{-}[d] & 2 \ar[dl]^1 &&
&& &&& 1 \ar[uu] \ar@{-}[d] & 1 && \\
\ar@{-}[rr] && 1 \ar@{-}[r] & 2 \ar[rrr] \ar[dl]^1 \ar@{-}[dddd] &&&
A
&& \ar@{-}[rrr] &&& *+[F]{2} \ar@{.>}[ddl] \ar[rrr] \ar@{-}[dddd] &&& A \\
&& 2 \ar[dl]^1 &&&&
&& && 1 &&&&\\
& 1 &&&& \HFKa &
&& && 1 \ar@{.>}[ddl] &&& \HFKm & \\
&&&&&&
&& &&&&&& \\
&&&&&&
&& &&&&&& \\
}
\]
\caption{
$\HFKa(10_{132})$ and $\HFKm(10_{132})$. Large numbers represent ranks
over $\Field$. The arrows in $\HFKa$ represent $\delta_1$ maps, with
ranks shown; the dotted arrows in $\HFKm$ represent multiplication by
$U$. The invariant $\tInvM$ lies in the boxed group.
}
\label{fig:10132HFK}
\end{figure}
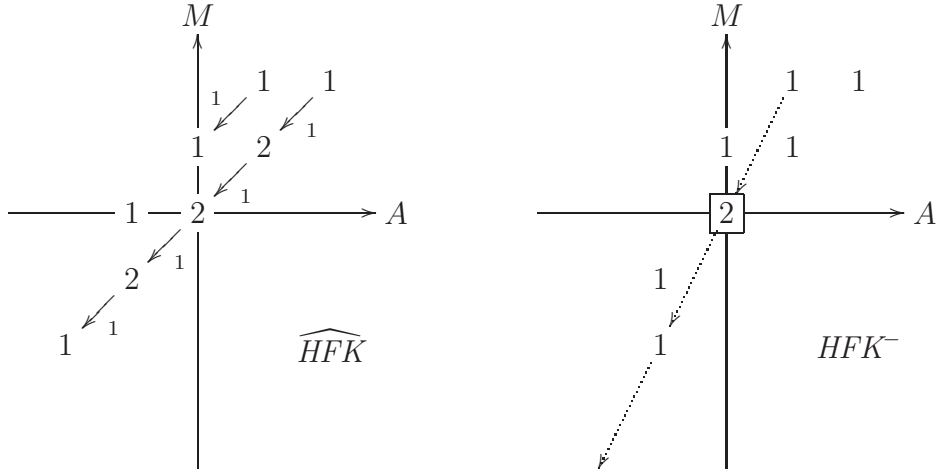

For the knot $K = m(10_{132})$, $\HFKa(m(K))$ and $\HFKm(m(K))$ are
plotted in Figure~\ref{fig:10132HFK}. Here $\HFKm_0(10_{132},0) =
\Field^2$, and indeed this is where $\tInvM$ sits for the two
transverse representatives of $m(10_{132})$.  The Floer homology
calculations are taken from Baldwin and Gillam's paper~\cite{BaldwinGillam}. More
precisely, they calculate $\HFKa$, together with its $\delta_1$ 
differential. In the given example, this can be used to determine $\HKm$, together
with its $U$ action, with the help of the following:

\begin{lemma}
  There is a spectral sequence starting at $H_*(\bigoplus_{t\leq -s}\HKa(K,t),\delta_1)$ and
  converging to $\HKm_{*-2s}(K,s)$.
  Moreover, there is a map of spectral sequences which induces the inclusion
  $$H_*(\bigoplus_{t\leq -s}\HKa_*(K,t),\delta_1)\longrightarrow H_*(\bigoplus_{t\leq -s+1}\HKa_*(K,t),\delta_1)$$
  on the $E_1$ page, and converges to the map
  $$U\colon \HKm_{*-2s}(K,s)\longrightarrow \HKm_{*-2s-2}(K,s-1).$$
\end{lemma}

\begin{proof}
  The spectral sequence comes from the following.  Consider
  the subcomplex $\Filta(K,s)$ of $\Ca(K)=\Cm(K)/(U_1=0)$ generated by
  those $\x$ with $A(\x)\leq s$. There are isomorphisms
  $$\phi_s\colon
  \HKm_{*}(K,s))\stackrel{\cong}{\longrightarrow}
  H_{*-2s}(\Filta_*(K,-s)),$$
  which fit into a commutative diagram
  \[
  \begin{CD}
    \HKm_{*}(K,s) @>{\phi_s}>> H_{*-2s}(\Filta_{*}(K,-s)) \\
    @V{U}VV   @VV{i}V \\
    \HKm_{*-2}(K,s-1)@>{\phi_{s-1}}>>H_{*-2s}(\Filta_{*}(K,-s+1)).
  \end{CD}
  \]
  This result is immediate for the holomorphic curves definition of the invariant~\cite{OS},
  see also the proof of Lemma A.2~\cite{OST}, for the result proved in the combinatorial context.
  
  We now consider the filtration of $\Filta(K,-s)$ by subcomplexes
  $\Filta(K,t)\subset \Filta(K,-s)$ with $t\leq -s$. The homology of the
  associated graded object is
  $\bigoplus_{t\leq -s}\HKa(K,-s)$, endowed with differential $\delta_1$.
\end{proof}

Glancing at the $\delta_1$ differential displayed on the
left-hand-side in Figure~\ref{fig:10132HFK}, we see at once
that the above spectral sequence collapses at the $E_2$ stage, and
hence that $\HKm$ is as shown on the right-hand-side of the same figure.

As for finding different transverse representatives in a candidate
knot class such as $m(10_{132})$, we found the program Gridlink
\cite{Gridlink}, and its ability to quickly produce Legendrian forms
of any reasonably small knot, to be very useful. In the case of
$m(10_{132})$, which is reversible, the Legendrian knot produced by
Gridlink has $\r=0$ and thus its orientation reverse gives a
transverse knot with the same $\tb$ and $\r$; these two knots are
the ones with different positive transverse pushoffs. If the trick of
reversing the orientation on a Legendrian knot with $\r=0$ does not
work, there are other ways to find candidates for different transverse knots.
For instance,
one can look in the front for a negative half-twist on two strands,
with one strand consisting of a downward-oriented zigzag, for which
the crossing is positive; one can then move the zigzag to the other
strand and produce another Legendrian knot with possibly different
transverse pushoff. See Figure~\ref{fig:local-grid-moves} or the dashed boxes in Figure~\ref{menasco} for an illustration (there are similar regions in
Figures~\ref{pretzel433} or~\ref{pretzel633}).

Finding possible candidate knot types for transverse knots which are
distinguished by whether or not $\delta_1 (\tInvA) = 0$ is
similar. Here we again need a knot $K$ such that some
$\HFKm_d(m(K),s)$ with $s=2d$ has rank at least $2$. In addition,
$\HFKa_d(m(K),s)$ should have rank at least $2$, and
$\HFKm_{d-1}(m(K),s-1)$ and $\HFKa_{d-1}(m(K),s-1)$ should be nonzero.

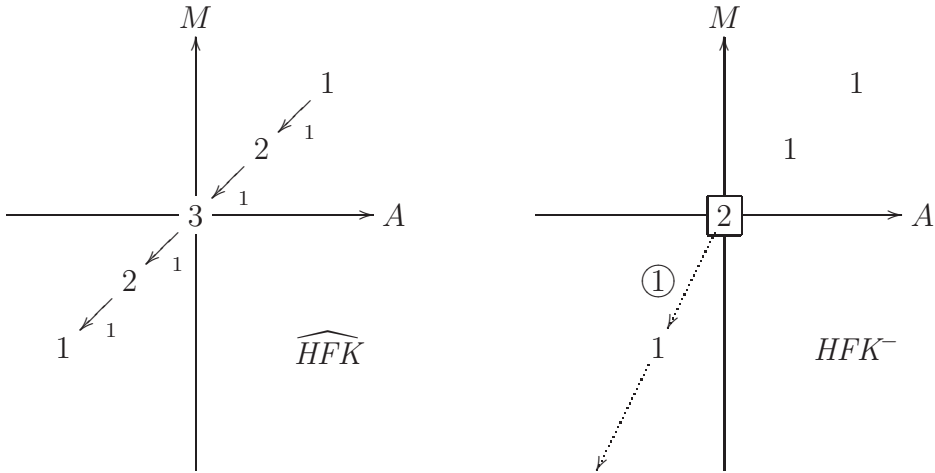
\begin{figure}
\xymatrix@=25pt@!0{ &&&M&&&
&& &&&M&&& \\
&&&&& 1 \ar[dl]^1 &
&& &&&&& 1 &\\
&&&& 2 \ar[dl]^1 &&
&& &&&& 1 && \\
\ar@{-}[rrr] &&& 3 \ar[uuu] \ar[rrr] \ar[dl]^1 \ar@{-}[dddd] &&& A
&& \ar@{-}[rrr] &&& *+[F]{2} \ar[uuu] \ar@{.>}[ddl] \ar[rrr] \ar@{-}[dddd] &&& A \\
&& 2 \ar[dl]^1 &&&&
&& && *+[o][F-]{1} &&&&\\
& 1 &&&& \HFKa&
&& && 1 \ar@{.>}[ddl] &&& \HFKm & \\
&&&&&&
&& &&&&&& \\
&&&&&&
&& &&&&&& \\
}
\caption{
$\HFKa(P(-4,-3,3))$ and $\HFKm(P(-4,-3,3))$. The invariant $\tInvM$
lies in the boxed group, $\delta_1\circ\tInvM$ in the circled group.
}
\label{fig:pretzelHFK}
\end{figure}

For $K = P(-4,-3,3) = m(10_{140})$, $\HFKa(m(K))$ and $\HFKm(m(K))$
are plotted in Figure~\ref{fig:pretzelHFK}. Note that
$\HFKm_0(m(K),0)$ has rank $2$; this is where $\tInvM$ sits for
the two transverse knots.

%% file: alg.tex
\section{The algorithm}
\label{sec:Algorithm}

We have seen in the last section applications of the transverse knot
invariant $\tInvA(L)\in\HFKa(\Mirror(L))$ (and also its image under $\delta_1$).

In practice, it is preferable to work with the finitely generated
chain complex $\CKaa(\Mirror(L);\Field)$ from Equation~\eqref{def:CKaa}. As in
Theorem~\ref{thm:KnotFloerHomology}, the homology groups of this
complex can be used to reconstruct $\HFKa(\Mirror(L))$ which, although it does
contain less information than $\HFKm(\Mirror(L))$, is sufficient for the
purposes of this article.

The strategy from~\cite{BaldwinGillam} for calculating knot Floer
homology can be adapted to the task at hand, trying to determine
whether or not $\tInvA$ is trivial. First observe that the question of
whether or not $\tInvA$ is homologically trivial is equivalent to the
question of whether or not the image $j_*(\tInvA)\in
\HKaa(\Mirror(L))$ is trivial, where here $j_*$ is the map on homology
induced by the projection $j\colon \CKa(\Mirror(L))\longrightarrow
\CKaa(\Mirror(L))$, cf.\ Proposition~\ref{prop:SpectralSequences} (we
are using here the statement that $j_1$ is an injection). Of course,
$j_*(\tInvA)$ is the homology class represented by the cycle $\xUR$,
thought of now as a homology class in $H_*(\CKaa(\Mirror(L)))$.

The limiting factor in determining knot Floer homology at the moment
is memory: for a knot with grid number $n$, one needs to keep track
of $n!$ generators. However, determining whether or not a given
cycle such as $\xUR$ is null-homologous requires less work than
calculating the ranks of all the homology groups, as much of the
chain complex is irrelevant to this problem. Indeed, this can be
formalized in the following algorithm.

Construct two sets
\[
  {\mathfrak A}=\coprod_{i=0}^\infty {\mathfrak A}_k \qquad\text{and}\qquad
  {\mathfrak B}=\coprod_{i=0}^\infty {\mathfrak B}_k
\]
defined inductively, as follows.  Let ${\mathfrak A}_0=\emptyset$ and
${\mathfrak B}_0=\{\xUR\}$.  Having built ${\mathfrak A}_k$
and ${\mathfrak B}_k$, we construct ${\mathfrak A}_{k+1}$ and
${\mathfrak B}_{k+1}$ as follows. Let ${\mathfrak A}_{k+1}$ consist of
all $\x\in\Gen(G)-{\mathfrak A}_{k}$ for which there is some $\y\in
{\mathfrak B}_{k}$ which appears with nonzero multiplicity in
$\daa(\x)$. Let ${\mathfrak B}_{k+1}$ consist of all
$\y\in\Gen(G)-{\mathfrak B}_k$ for which there is some $\x\in
{\mathfrak A}_{k+1}$ for which $\y$ appears with nonzero multiplicity
in $\daa(\x)$.  The important point here is that the construction of
${\mathfrak A}_{k+1}$ and ${\mathfrak B}_{k+1}$ requires keeping
track of the grid states representing only ${\mathfrak A}_k$ and ${\mathfrak
  B}_k$ (rather than all the ${\mathfrak A}_i$ and ${\mathfrak B}_i$
for $i=1,\dots,k$).

Let $A$ (respectively $B$, $A_k$, $B_k$) be the free vector space over $\Field$ generated by elements in
${\mathfrak A}$ (respectively ${\mathfrak B}$, ${\mathfrak A}_k$, ${\mathfrak B}_k$), and form the chain complex $C'=A\oplus
B$, endowed with the differential
$$D\colon A \longrightarrow B$$ gotten by counting rectangles as in the
definition of $\daa$.
By construction, there is a natural quotient map
$$Q\colon \CKaa\longrightarrow C'.$$
It is easy to see that $\xUR$ is homologically trivial
if and only if $Q(\xUR)$ is; so, since $C'$ is significantly smaller than
$\CKaa$ (in particular, it generated by elements with Maslov grading
$d$ and $d+1$, where $d=\sl(L)+1$), we work with it instead.

To determine whether or not the given element $\xUR$ is nontrivial in
$C'$, we proceed as follows. First, we enlarge $C'$ to a chain complex
$C''=A'\oplus B$ with one additional (distinguished) generator $a_0\in
A'$ (i.e. $A'=A\oplus\Field$) with $D'(a_0)=\xUR$. (We can think of
$\xUR$ as inducing a chain map from $\Field$ to $C'$, and hence that
$C''$ is the mapping cone of this map.)  As in Baldwin and
Gillam~\cite{BaldwinGillam}, we view the differential on the complex
$C''$ as giving a graph on the generating set of $C''$, drawing an
edge from a generator $a$ of $C''$ to $b$ in $C''$ whenever $b$
appears with nonzero multiplicity in $D'(a)$. As in their scheme,
given an edge $e$ from $a$ to $b$, we can ``contract'' it (and
reduce the number of generators of our complex by two) without
affecting the homology, as follows. Draw additional edges from
$a'$ to $b'$ (and then cancel identical edges in pairs), for all 
$a'$ for which there is an edge from $a'$ to $b$ and all $b'$ for which there is an edge from $a$ to $b'$.
This has the effect of a change of basis $a' \mapsto a' + a$, $b' \mapsto b' + b$ for all such $a'$, $b'$. After this change of basis, the only edge involving $a$ or $b$ is the edge between them, and we can then contract the complex by deleting $a$, $b$, and the edge between them.

It will be
important for us to perform these contractions in a controlled manner.
Let $C^0$ be the initial complex $C''$. We construct now a finite
sequence of contractions to obtain a finite sequence of complexes
$\{C^k\}_{k=0}^m$ as follows. Given $C^k$, we consider the
distinguished element $a_0$. If there are no edges leaving $a_0$, then
our sequence terminates, and we conclude that $\xUR$ had to be
homologically trivial in $C$, and hence also in $\CKaa$. If there are
edges leaving $a_0$, we ask if there are edges $e$ connecting some
$a\neq a_0$ to some $b$ which is also the endpoint of a different edge
leaving $a_0$. If there is no such edge $e$, then our sequence terminates,
and we conclude that the original element $\xUR$ had to be
homologically nontrivial in $C$, and hence also in $\CKaa$. Otherwise,
we let $C^{k+1}$ be the complex obtained from $C$ after contracting~$e$.

Finally, we remark that some time is saved if we perform both
operations simultaneously: before building the next level of the complex
$A_{k+1}$ and $B_{k+1}$ from $A_k$ and $B_k$, we contract all possible
edges which point from $A_{k+1}$ into $B_{k}$. This is, in fact, the algorithm we implemented for performing the calculations from the present paper.

We indicate how this works in a particular example, the pretzel knot
$P(-4,-3,3)$ represented by
\begin{align*}
X = (9,8,1,4,6,5,7,2,3), ~O = (4,2,5,7,9,8,3,6,1)
\end{align*}
investigated in the last section. The computation is slightly
easier, and the algorithm equally well demonstrated, if we
investigate $\xLL$ rather than $\xUR$. In this case, it turns out
that the above algorithm encounters only $8$ grid states (of
the $9!=362880$ grid states which generate $\CKaa(G)$).

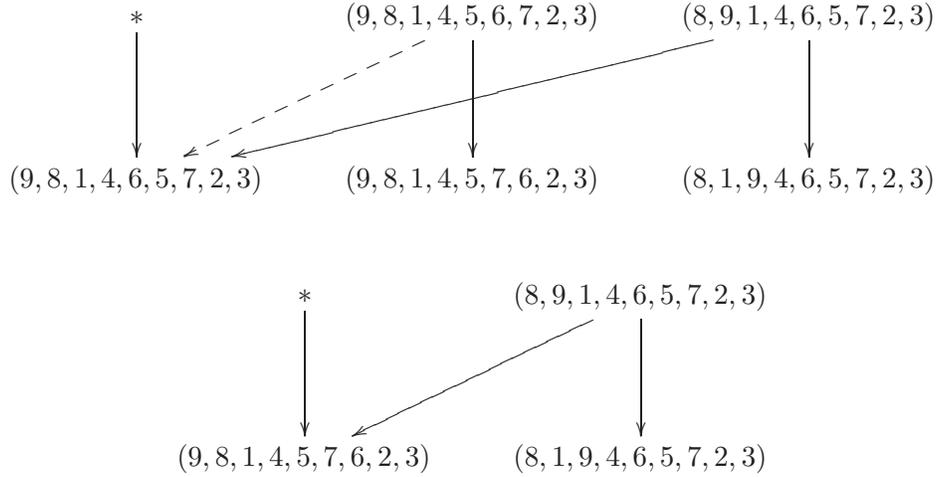
\begin{figure}[ht]
\[\small
\xymatrix@R+20pt{\ast\ar[d] & (9,8,1,4,5,6,7,2,3)\ar[d]\ar@{-->}[dl] & (8,9,1,4,6,5,7,2,3)\ar[d]\ar[dll]\\
(9,8,1,4,6,5,7,2,3) & (9,8,1,4,5,7,6,2,3) & (8,1,9,4,6,5,7,2,3)}
\]\vspace{20pt}
\[\small
\xymatrix@R+20pt{\ast\ar[d] & (8,9,1,4,6,5,7,2,3)\ar[d]\ar[dl]\\
(9,8,1,4,5,7,6,2,3) & (8,1,9,4,6,5,7,2,3)}
\]
\caption{\label{fig:Complex}
{\bf{Contracting edges in the chain complex for $P(-4,-3,3)$.}}
At the first stage of building the complex for $P(-4,-3,3)$, we obtain the complex
pictured on the top. Here, the asterisk denotes the distinguished element $a_0$.
Contracting the dotted edge, we obtain the second complex, shown on the bottom.}
\end{figure}

We start with the initial state $\xLL$, which we could alternatively
denote by its $y$-coordinates $(9,8,1,4,6,5,7,2,3)$. Thus, we start
with the chain complex $A_0'$ generated by the distinguished element
$a_0$, and $B_0$, generated by
$${\mathfrak B}_0=\{(9,8,1,4,6,5,7,2,3)\}.$$
A casual glance at the
grid picture reveals exactly two rectangles pointing into this
grid state, and these are rectangles leaving states
$$ {\mathfrak A}_1=
\{(9,8,1,4,5,6,7,2,3),(8,9,1,4,6,5,7,2,3)\}. $$ For each of the
above states $a\in {\mathfrak A}_1$, there are exactly two empty
rectangles leaving $a$; one goes back to $\xLL$, and the other points
to a new state in
$${\mathfrak B}_1=\{(9,8,1,4,5,7,6,2,3),(8,1,9,4,6,5,7,2,3)\}.$$
Thus, so far, we have the complex pictured on the top in
Figure~\ref{fig:Complex}.

\begin{figure}[ht]
\mbox{\vbox{\epsfbox{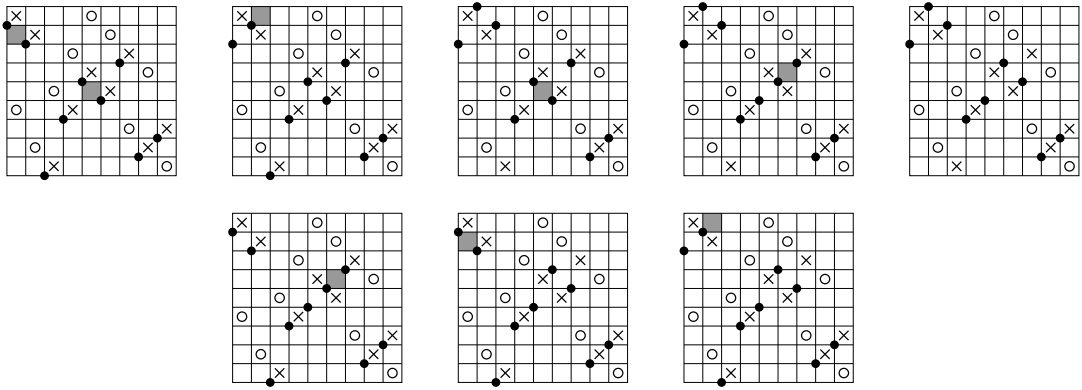}}}
\caption{\label{fig:PretzelGridStates}
  {\bf{Finding states for the chain complex for $P(-4,-3,3)$.}}
  The columns represent states in ${\mathfrak B}_0$, ${\mathfrak A}_1$,
  ${\mathfrak B}_1$, ${\mathfrak A}_2$, ${\mathfrak B}_2$ respectively.
  Shaded rectangles connect states in a given column to states in the
  next column.}
\end{figure}

Contracting one edge pointing back to the generator
$(9,8,1,4,6,5,7,2,3)\in {\mathfrak B}_0$, we obtain a smaller complex,
having thrown out the generator of ${\mathfrak B}_0$. But we have no
further need for such a generator: in building ${\mathfrak A}_2$, we
need only remember the states in ${\mathfrak B}_1$. Proceeding in this
manner, the persistent reader can verify that
\begin{align*}
{\mathfrak A}_2&=
\left\{
(8,9,1,4,5,7,6,2,3),(8,1,9,4,5,6,7,2,3)
\right\} \\
{\mathfrak B}_2 &=
\{(8,1,9,4,5,7,6,2,3)\}.
\end{align*}
There are no new states (i.e. states not already in ${\mathfrak
A}_2$) pointing into ${\mathfrak B}_2$; thus we need only calculate
the homology of what we have so far, to determine that $a_0$ is not
a cycle, and hence that $\xLL$ was homologically nontrivial.

A mild modification of the above procedure applies when determining whether
 $\delta_1 (\tInvA)$ is null-homologous. In this case, we start
with ${\mathfrak A}_0 = \emptyset$ and 
${\mathfrak B}_0$ consisting of all of the terms in $\daa_1(\xUR)$,
that is, containing exactly one $X$:
\[
{\mathfrak B}_0 = \sum \{\y \mid \textstyle\y \in \Gen, \exists !~ r \in \EmptyRect(\x,\y) ~\text{with}~
   \sum O_i(r) = 0, \sum X_i(r) = 1 \}.
\]
\noindent
We build the complex $C'$ from here as before. We enlarge this to $C''$ by adding
one additional distinguished generator $a_0$ whose differential consists of all the terms
in ${\mathfrak B}_0$. Arguing as before, we have that
$\delta_1(\tInvA)$ is trivial if and only if $a_0$ has no edges pointing out of it (after
all other edges have been contracted).

Although the pretzel example we illustrated above had a reasonably
small complex, most of the other examples in this paper are quite
involved; and hence, we implemented the above algorithm in a C
program.

%% file: nat.tex
\section{Naturality questions}
\label{sec:naturality}

We end with a conjecture that would, if true, greatly increase the
power of the transverse invariant $\tInv$.  Consider a sequence~$S$ of
elementary grid diagram moves, including symmetries of the grid
diagram that fix the vertical axis,
that start and end at the same grid diagram~$G$ for a knot~$K$.  This
sequence of grid diagrams gives a path in the space of embeddings
of~$K$, and thus an element $[S]$ in $\Mod^+(S^3,K)$, the mapping
class group of~$S^3$ relative to~$K$, preserving the orientations of
both $S^3$ and~$K$.

\begin{conjecture}\label{conj:natural}
  The map $\Phi(S) : \HFKm(G) \rightarrow \HFKm(G)$ induced by~$S$ on
  Heegaard-Floer homology depends only on $[S] \in \Mod^+(S^3,K)$, up
  to a possible overall sign.
\end{conjecture}

More generally, one might conjecture that any cobordism between knots
gives a well-defined map on the homology:

\begin{conjecture}
  For any pair of oriented knots $K, K' \subset S^3$ and oriented
  surface $\Sigma \subset S^3 \times [0,1]$ so that $\partial\Sigma =
  (-K \times \{0\}) \cup (K' \times \{1\})$ (where $-K$ is $K$ with
  the orientation reversed), there is a map
  \[
  \Phi(\Sigma) : \HFKm(K) \rightarrow \HFKm(K')
  \]
  which depends only on the isotopy class of $\Sigma$.
\end{conjecture}

(The analogous theorem is true for Khovanov
homology~\cite{Jacobsson,MW:Functoriality}.)

\begin{figure}
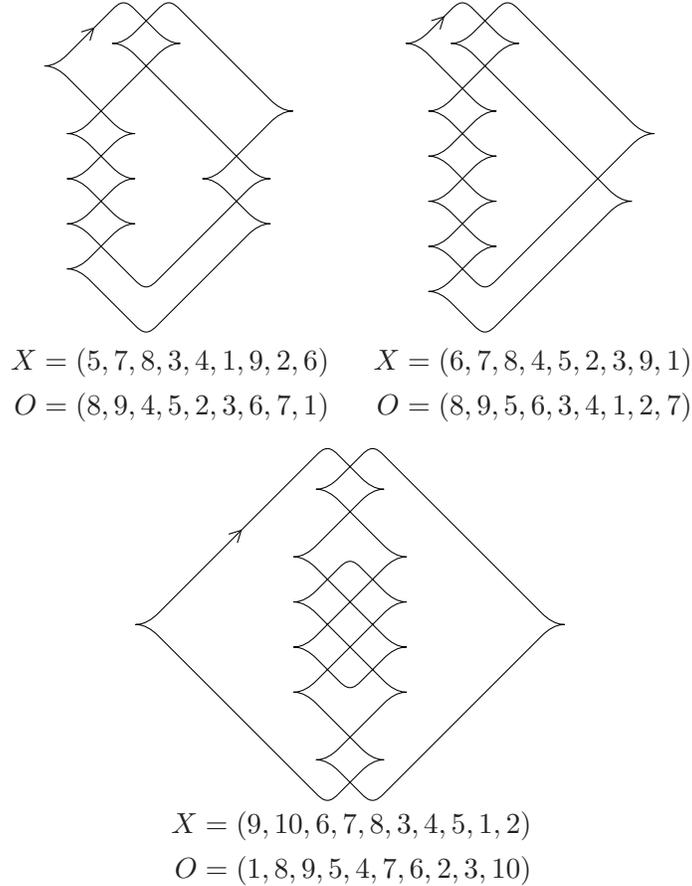

  \[
  \TransDiag{draws/knots.42}
    {5,7,8,3,4,1,9,2,6}
    {8,9,4,5,2,3,6,7,1}
  \quad
  \TransDiag{draws/knots.40}
    {6,7,8,4,5,2,3,9,1}
    {8,9,5,6,3,4,1,2,7}
  \]
  \[
  \mathcenter{\TransDiag{draws/knots.46}
    {9,10,6,7,8,3,4,5,1,2}
    {1,8,9,5,4,7,6,2,3,10}}
  \]
  \caption{Legendrian knots from grid diagrams $G_1$ for $E(2,4)$
    (top left) and $G_2$ for $E(1,5)$ (top
    right) of type $7_2$.  An alternate, symmetric form $G_3$ of
    $E(1,5)$ is shown below.}
  \label{fig:7_2}
\end{figure}

If Conjecture~\ref{conj:natural} were true, the transverse invariant
$\tInvM(T)$ of a transverse knot~$T$ of topological type~$K$ would
be well-defined as an element of $\HFKm(K)/\Mod^+(S^3,K)$.  Since
the mapping class groups $\Mod^+(S^3,K)$ are generally smaller than
the group of all automorphisms of $\HFKm(K)$, the invariant would
become stronger.  For instance, let $E(1,5)$ and $E(2,4)$ denote two
of the ``Eliashberg knots'' considered by Epstein, Fuchs, and Meyer
\cite{EpsteinFuchsMeyer}\footnote{Warning: \cite{EpsteinFuchsMeyer}
  uses the opposite convention for transverse pushoffs; their
  $L^\pm$ is our $L^\mp$.} and
shown in Figure~\ref{fig:7_2}.
Then we have the following.

\begin{proposition}
  If Conjecture~\ref{conj:natural} is true, then $E(1,5)^+$ and
  $E(2,4)^+$ are not transversely isotopic.
\label{prop:naturality}
\end{proposition}

\begin{figure}
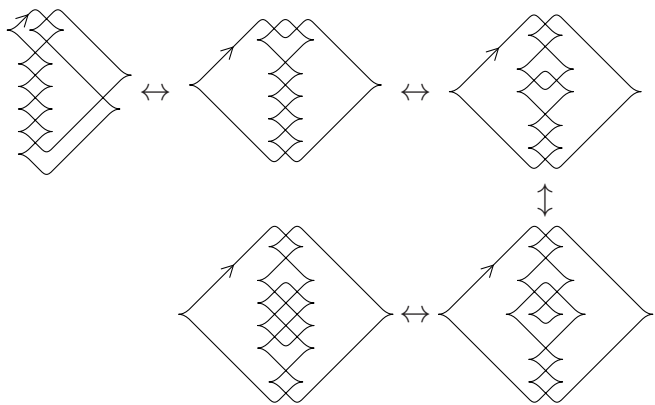

  \[
  \begin{tabular}{@{}c@{}c@{}c@{}c@{}c@{}}
    $\mfigb{draws/knots.80}$&$\leftrightarrow$&$\mfigb{draws/knots.81}$&$\leftrightarrow$&$\mfigb{draws/knots.82}$\\
    &&&&$\updownarrow$\\
    &&$\mfigb{draws/knots.84}$&$\leftrightarrow$&$\mfigb{draws/knots.83}$
  \end{tabular}
  \]
  \caption{A Legendrian isotopy between $G_2$ and $G_3$, the two forms
    of $E(5,1)$.  Starting from the upper left and moving clockwise,
    rotate (in the toroidal grid diagram) to the right one step, then rotate to the right and down two steps each,
    then perform a stabilization of type $\ONW$, and finally perform
    two commutations.}
  \label{fig:isotopy-G2-G3}
\end{figure}

\begin{figure}
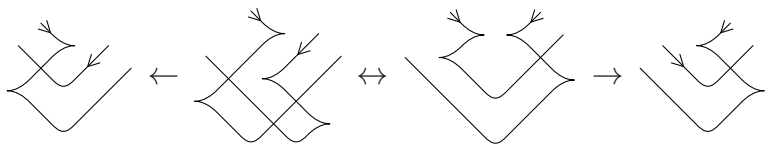

  \[
  \mfigb{draws/knots.70}\leftarrow\mfigb{draws/knots.71}\leftrightarrow\mfigb{draws/knots.72}\rightarrow\mfigb{draws/knots.73}
  \]
  \caption{A local sequence of grid moves that goes from $G_1$ to
  $G_2$. In the Legendrian category, these comprise a positive
  stabilization and destabilization. This same sequence can also be
  used to
  relate the Legendrian pretzel knot pairs and the Etnyre-Honda pair
  considered in Section~\ref{sec:Examples}.
}
  \label{fig:local-grid-moves}
\end{figure}

\begin{proof}
  Both knots are of topological type $7_2$.  The mapping class group
  $\Mod^+(S^3,7_2)$ is $\Z/2\Z$~\cite{RW}; the nontrivial element~$\phi$ is
  one that exists for every two-bridge knot: if the knot is in bridge
  position with respect to a sphere~$S$, $\phi$ interchanges the two
  positive and two negative intersections with~$S$ (so preserving the
  orientation of the knot).
  
  Let $G_1 = E(2,4)$, $G_2 = E(1,5)$, $G_3$ be the diagrams shown in Figure~\ref{fig:7_2}, and note that one
  representative for $\phi$ is rotation by $180^\circ$ on $G_3$.
  This symmetry interchanges $[\z^+(G_3)]$ and $[\z^-(G_3)]$, which
  therefore form an orbit of $\Mod^+(S^3,7_2)$.  Now $G_2$ and $G_3$ are Legendrian isotopic by Figure~\ref{fig:isotopy-G2-G3}, and $[\z^\pm]$ are preserved by Legendrian isotopy; hence $[\z^+(G_2)]$ and $[\z^-(G_2)]$ form an orbit of $\Mod^+(S^3,7_2)$.
  
  On the other hand,
  $G_1$ and $G_2$ can be related by a sequence of grid moves, as shown in Figure~\ref{fig:local-grid-moves}. It is straightforward to check by hand, using the quasi-isomorphisms under grid moves from \cite{MOST}, that this sequence takes $[\z^+(G_1)]$ to $(2,7,8,4,5,9,3,6,1) + (4,7,8,9,5,2,3,6,1)$ in the complex $\CKaa$ for $G_2$ (in the notation from Section~\ref{sec:Algorithm}). By inspection, this sum is not homologous in $\CKaa$ to $[\z^+(G_2)]$ or $[\z^-(G_2)]$. On the assumption of
  Conjecture~\ref{conj:natural}, it follows that the positive
  transverse pushoffs $E(1,5)^+$ and $E(2,4)^+$ are distinct.
\end{proof}

\noindent By comparison, note that, by
\cite[Theorem~2.2]{EpsteinFuchsMeyer} (or see
Figure~\ref{fig:local-grid-moves}), $E(1,5)^-$ and $E(2,4)^-$ are
transversely isotopic.

We close this section with some speculation about transverse twist
knots. The knots $E(k,l)$ of \cite{EpsteinFuchsMeyer}, which are
twist knots of crossing number $n+1$ if $k+l=n$, generalize
Chekanov's celebrated examples of nonisotopic Legendrian $5_2$
knots \cite{Chekanov}. More specifically, $E(k,l)$ and $E(k',l')$ are Legendrian isotopic if and only if $\{k,l\} = \{k',l'\}$ \cite{EpsteinFuchsMeyer}, and the $E(k,l)$
provide candidates for nonisotopic transverse knots. It
is conjectured in \cite{EpsteinFuchsMeyer} that $E(k,l)^+$ and
$E(k-1,l+1)^+$ are not transversely isotopic in general whenever $l$
is even. This turns out not to be true: for instance, the transverse
$5_2$ knots $E(2,2)^+$ and $E(1,3)^+$ are isotopic.

In general, one can show that if $n$ is odd, then the oriented knots
$$E(1,n-1), E(2,n-2), \dots, E(n-1,1)$$ are all Legendrian isotopic
after one negative stabilization (and also all isotopic after one
positive stabilization). Analogously, if $k,l\geq 2$ are both even,
then $E(k,l)$, $E(k+1,l-1)$, $E(l-1,k+1)$, and $E(l-2,k+2)$ (this
last assuming $l \geq 4$) are all Legendrian isotopic after one
negative stabilization. Even with this in mind, however, it is still
possible that if $n$ is even, several of $E(1,n-1), E(2,n-2), \dots,
E(n-1,1)$ are not negatively stably isotopic.

\begin{conjecture}
If $n$ is even, the $\lceil n/4 \rceil$ transverse knots
\[
E(1,n-1)^+, E(3,n-3)^+, \dots, E(2\lceil n/2 \rceil - 1, 2\lfloor
n/2 \rfloor + 1)^+,
\]
all with $\sl=1$, are pairwise transversely nonisotopic.
\label{conj:twist}
\end{conjecture}

We note that the Legendrian knots $E(1,n-1), E(3,n-3), \dots,
E(2\lceil n/2 \rceil - 1, 2\lfloor n/2 \rfloor + 1)$ are hardly
arbitrarily chosen: they are particularly simple examples of the
Legendrian satellite construction \cite{NgTraynor}. More precisely,
they are Whitehead doubles of the unknot obtained by taking the
Legendrian satellites of unknots with Thurston-Bennequin number
$-n/2$ and the Legendrian Whitehead knot $W_0$ described in the
appendix of \cite{NgTraynor}. It is straightforward to check that
these Legendrian Whitehead doubles are unchanged up to isotopy by
reversing the orientation of the underlying unknot; note that the
number of different unoriented Legendrian unknots with $\tb=-n/2$ is
$\lceil n/4 \rceil$ \cite{EliashbergFraser}.

Considering the Alexander polynomial, and using the structure of knot
Floer homology of two-bridge knots~\cite{RasmussenTwoBridge}, see
also~\cite{AltKnots}, one sees that the transverse invariant $\tInv$
for the knots in Conjecture~\ref{conj:twist} lies in an $\HFKa$ group
of rank $n/2$. Just as for $7_2$, the mapping class group of the
underlying topological twist knot is $\Z/2\Z$; quotienting by the
$\Z/2\Z$ action yields a group of rank $\lceil n/4 \rceil$. It does
not seem unreasonable to guess that each of the $\lceil n/4 \rceil$
transverse knots from Conjecture~\ref{conj:twist} maps under $\tInv$
to a different generator in the quotient.


%% file: paper.bbl
\begin{thebibliography}{10}

\bibitem{BaldwinGillam}
J.~A. Baldwin and W.~D. Gillam.
\newblock Computations of {H}eegaard-{F}loer knot homology.
\newblock math.GT/0610167.

\bibitem{BirmanMenascoNote}
J.~S. Birman and W.~M. Menasco.
\newblock A note on closed 3-braids.
\newblock arXiv:0802.1072.

\bibitem{BirmanMenasco}
J.~S. Birman and W.~M. Menasco.
\newblock Stabilization in the braid groups. {II}. {T}ransversal simplicity of
  knots.
\newblock {\em Geom. Topol.}, 10:1425--1452 (electronic), 2006.
\newblock math.GT/0310280.

\bibitem{Brunn}
H.~Brunn.
\newblock {\"U}ber verknotete {K}urven.
\newblock {\em Verhandlungen des Internationalen Math. Kongresses (Z{\"u}rich
  1897)}, pages 256--259, 1898.

\bibitem{Chekanov}
Y.~Chekanov.
\newblock Differential algebra of {L}egendrian links.
\newblock {\em Invent. Math.}, 150(3):441--483, 2002.

\bibitem{MW:Functoriality}
D.\ Clark, S.\ Morrison, and K.\ Walker.
\newblock Fixing the functoriality of {K}hovanov homology.
\newblock math.GT/0701339.

\bibitem{Cromwell}
P.~R. Cromwell.
\newblock Embedding knots and links in an open book. {I}. {B}asic properties.
\newblock {\em Topology Appl.}, 64(1):37--58, 1995.

\bibitem{Gridlink}
M.~Culler.
\newblock Gridlink, 2006--2007.
\newblock Available from \url{http://www.math.uic.edu/~culler/gridlink/}.

\bibitem{Dynnikov}
I.~A. Dynnikov.
\newblock Arc-presentations of links: monotonic simplification.
\newblock {\em Fund. Math.}, 190:29--76, 2006.
\newblock math.GT/0208153.

\bibitem{EliashbergFraser}
Y.~Eliashberg and M.~Fraser.
\newblock Classification of topologically trivial {L}egendrian knots.
\newblock In {\em Geometry, topology, and dynamics (Montreal, PQ, 1995)},
  volume~15 of {\em CRM Proc. Lecture Notes}, pages 17--51. Amer. Math. Soc.,
  Providence, RI, 1998.

\bibitem{EpsteinFuchsMeyer}
J.~Epstein, D.~Fuchs, and M.~Meyer.
\newblock Chekanov-{E}liashberg invariants and transverse approximations of
  {L}egendrian knots.
\newblock {\em Pacific J. Math.}, 201(1):89--106, 2001.

\bibitem{Etnyresurvey}
J.~B. Etnyre.
\newblock Legendrian and transversal knots.
\newblock In {\em Handbook of knot theory}, pages 105--185. Elsevier B. V.,
  Amsterdam, 2005.
\newblock math.SG/0306256.

\bibitem{EtnyreHondaJSG}
J.~B. Etnyre and K.~Honda.
\newblock Knots and contact geometry. {I}. {T}orus knots and the figure eight
  knot.
\newblock {\em J. Symplectic Geom.}, 1(1):63--120, 2001.
\newblock math.GT/0006112.

\bibitem{EtnyreHonda}
J.~B. Etnyre and K.~Honda.
\newblock Cabling and transverse simplicity.
\newblock {\em Ann. of Math. (2)}, 162(3):1305--1333, 2005.
\newblock math.SG/0306330.

\bibitem{ENS}
J.~B. Etnyre, L.~L. Ng, and J.~M. Sabloff.
\newblock Invariants of {L}egendrian knots and coherent orientations.
\newblock {\em J. Symplectic Geom.}, 1(2):321--367, 2002.
\newblock math.SG/0101145.

\bibitem{Jacobsson}
M.~Jacobsson.
\newblock An invariant of link cobordisms from {K}hovanov homology.
\newblock {\em Algebr. Geom. Topol.}, 4:1211--1251, 2004.
\newblock math.GT/0206303.

\bibitem{MOS}
C.~Manolescu, P.~S. Ozsv{\'a}th, and S.~Sarkar.
\newblock A combinatorial description of knot {F}loer homology.
\newblock math.GT/0607691. To appear in {\em Ann. Math.}

\bibitem{MOST}
C.~Manolescu, P.~S. Ozsv{\'a}th, Z.~Szab{\'o}, and D.~P. Thurston.
\newblock On combinatorial link {F}loer homology.
\newblock {\em Geom. Top.}, 11:2339--2412, 2007.
\newblock arXiv:math/0610559.

\bibitem{UsersGuide}
J.~McCleary.
\newblock {\em A user's guide to spectral sequences}, volume~58 of {\em
  Cambridge Studies in Advanced Mathematics}.
\newblock Cambridge University Press, Cambridge, second edition, 2001.

\bibitem{MenascoMatsuda}
W.~W. Menasco and H.~Matsuda.
\newblock An addendum on iterated torus knots (appendix).
\newblock math.GT/0610566.

\bibitem{NgKhovanov}
L.~Ng.
\newblock A {L}egendrian {T}hurston-{B}ennequin bound from {K}hovanov homology.
\newblock {\em Algebr. Geom. Topol.}, 5:1637--1653 (electronic), 2005.
\newblock math.GT/0508649.

\bibitem{NgTraynor}
L.~Ng and L.~Traynor.
\newblock Legendrian solid-torus links.
\newblock {\em J. Symplectic Geom.}, 2(3):411--443, 2004.
\newblock math.SG/0407068.

\bibitem{NgCLI}
L.~L. Ng.
\newblock Computable {L}egendrian invariants.
\newblock {\em Topology}, 42(1):55--82, 2003.
\newblock math.GT/0011265.

\bibitem{OS}
P.~Ozsv{\'a}th and Z.~Szab{\'o}.
\newblock Holomorphic disks and knot invariants.
\newblock {\em Adv. Math.}, 186(1):58--116, 2004.
\newblock math.GT/0209056.

\bibitem{OST}
P.~Ozsv\'ath, Z.~Szab\'o, and D.~Thurston.
\newblock Legendrian knots, transverse knots and combinatorial {F}loer
  homology.
\newblock math.GT/0611841.

\bibitem{AltKnots}
P.~S. Ozsv{\'a}th and Z.~Szab{\'o}.
\newblock Heegaard {F}loer homology and alternating knots.
\newblock {\em Geom. Topol.}, 7:225--254 (electronic), 2003.
\newblock math.GT/0209149.

\bibitem{HolDisk}
P.~S. Ozsv{\'a}th and Z.~Szab{\'o}.
\newblock Holomorphic disks and topological invariants for closed
  three-manifolds.
\newblock {\em Ann. of Math. (2)}, 159(3):1027--1158, 2004.
\newblock math.SG/0101206.

\bibitem{Plam:HFK}
O.~Plamenevskaya.
\newblock Bounds for the {T}hurston-{B}ennequin number from {F}loer homology.
\newblock {\em Algebr. Geom. Topol.}, 4:399--406, 2004.
\newblock math.SG/0311090.

\bibitem{Plamenevskaya}
O.~Plamenevskaya.
\newblock Transverse knots and {K}hovanov homology.
\newblock {\em Math. Res. Lett.}, 13(4):571--586, 2006.
\newblock math.GT/0412184.

\bibitem{RasmussenTwoBridge}
J.~A. Rasmussen.
\newblock Floer homology of surgeries on two-bridge knots.
\newblock {\em Algebr. Geom. Topol.}, 2:757--789 (electronic), 2002.
\newblock math.GT/0204056.

\bibitem{Rasmussen}
J.~A. Rasmussen.
\newblock {\em Floer homology and knot complements}.
\newblock PhD thesis, Harvard University, 2003.
\newblock math.GT/0611841.

\bibitem{RW}
A.~Reid and G.~Walsh.
\newblock Commesurability classes of 2--bridge knot complements.
\newblock arXiv:math.GT/0612473.

\bibitem{SarkarWang}
S.~Sarkar and J.~Wang.
\newblock {A combinatorial description of some {H}eegaard {F}loer homologies}.
\newblock math.GT/0607777.

\end{thebibliography}
